\RequirePackage{amsthm}
\documentclass[sn-mathphys-num]{sn-jnl}

\usepackage{graphicx}%
\usepackage{multirow}%
\usepackage{amsmath,amssymb,amsfonts}%
\usepackage{amsthm}%
\usepackage{mathrsfs}%
\usepackage[title]{appendix}%
\usepackage{xcolor}%
\usepackage{textcomp}%
\usepackage{manyfoot}%
\usepackage{booktabs}%
\usepackage{algorithm}%
\usepackage{algorithmicx}%
\usepackage{algpseudocode}%
\usepackage{listings}%
 \usepackage{enumitem}

\newtheorem{theorem}{Theorem}[section]
\newtheorem{lemma}[theorem]{Lemma}

\newtheorem{proposition}[theorem]{Proposition}

\theoremstyle{definition}
\newtheorem{definition}[theorem]{Definition}
\newtheorem{example}[theorem]{Example}

\theoremstyle{remark}
\newtheorem{remark}[theorem]{Remark}

\numberwithin{equation}{section}

\providecommand{\abs}[1]{\lvert#1\rvert}
\providecommand{\norm}[1]{\lVert#1\rVert}

\begin{document}

\title[Cartan--Thullen theorem, Levi problem for generalised convexity]{Cartan--Thullen theorem and Levi problem in context of generalised convexity}

\author[1,2]{\fnm{Krzysztof J.} \sur{Ciosmak} \href{https://orcid.org/0000-0001-9571-1160}{\includegraphics{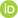}}}\email{k.ciosmak@utoronto.ca}

\affil[1]{\orgdiv{Department of Mathematics}, \orgname{University of Toronto}, \orgaddress{\street{Bahen Centre, 40 St. George St., Room 6290}, \city{Toronto}, \postcode{Ontario, M5S 2E4}, \country{Canada}}}
\affil[2]{\orgname{Fields Institute for Research in Mathematical Sciences}, \orgaddress{\street{222 College Street}, \city{Toronto}, \postcode{Ontario, M5T 3J1}, \country{Canada}}}


\abstract{We demonstrate that the Cartan--Thullen theorem and its generalisation to the context of generalised convexity, which we establish herein, can be regarded as consequences of the classical theorems of functional analysis: the Banach--Steinhaus theorem and the Banach--Alaoglu theorem.
Furthermore, we characterise the domains of holomorphy, and their generalisations, as the spaces that are complete, or as the spaces exhaustible by suitably defined polytopes.
We also provide an abstract analogue of the Levi problem and its elementary resolution. 
Our results allow for a novel characterisation of Stein spaces as the holomorphically complete spaces, as well as a proof that the Bremermann--Lelong lemma is equivalent to the positive answer to the Levi problem.
Another contribution of ours is the introduction of the analogues of the notions of the complex analysis to the setting of generalised convexity.
}



\maketitle

\bmhead{Acknowledgments}
I would like to express my thanks to Robert McCann for bringing to my attention references regarding generalised convexity. I thank also Thomas Bloom for the comments that helped improving the manuscript. Part of the research presented in this paper was performed when the author was a post-doctoral fellow at the University of Oxford, supported by the ERC Starting Grant CURVATURE, grant agreement No. 802689.

\maketitle

\section{Introduction}\label{s:intro}

Our first aim is to show that the Cartan--Thullen theorem is an instance of a novel theorem in abstract generalised convexity, whose proof relies solely on functional analytic tools, such as the Banach-Alaoglu theorem and the Banach--Steinhaus uniform boundedness principle. 

The second aim is to formulate the Levi problem in this abstract setting, and to resolve it.  Our resolution relies on elementary soft analytic arguments. The resolution implies several results presented in the book of H\"ormander \cite{Hormander2007} as particular instances.

We therefore show that the analogies between classical convexity and holomorphic convexity, as exposed in \cite{Bremermann19562} or in \cite{Range2012}, pertain to a much broader class of notions of generalised convexity.

The approaches to the original Levi problem either used  homological methods, due to Oka \cite{Oka1950}, Bremermann \cite{Bremermann1954} and Norguet \cite{Norguet1954}, or involved partial differential equations, due to H\"ormander \cite{Hormander1965}.
In our developments the exhaustion function, which is the matter of concern in the Levi problem, is directly constructed using the very definition of holomorphic convexity, or the analogue of such convexity in the generalised convexity setting.\footnote{In the general abstract setting of a topological space $\Omega$ and a linear space of continuous functions $\mathcal{A}$ on $\Omega$, studied in this paper, there is no clear notion of functions in $\mathcal{A}$ on the open subsets of $\Omega$ -- unless  the restrictions are considered, in which case the Levi problem becomes trivial. Therefore we consider the formulation concerning the exhaustion functions as the most appropriate generalisation of the Levi problem to our abstract setting. However, we do not obtain a new resolution of the original Levi problem.}

A further contribution presented in this paper is the introduction of the analogues of the notions of domain of holomorphy and domain of existence to the setting of generalised convexity. Also notions of completeness and exhaustibility by polytopes are provided and their equivalence to the completeness with respect to the considered class of functions is shown. In the complex analytic setting it is proven that the holomorphic completeness is equivalent to Steinness of  a complex analytic space. This shows that the intention of Stein to name the Stein spaces as holomorphically complete, see \cite{Stein1951}, was accurate. 

\subsection{The Cartan--Thullen theorem and the Levi problem}

The classical result of Cartan and Thullen is concerned with  holomorphically convex subsets of $\mathbb{C}^n$. 

\begin{definition}
Let $\Omega\subset\mathbb{C}^n$ be an open set. For $B\subset\Omega$ its closed holomorphically convex hull is defined as 
\begin{equation*}
\{z\in \Omega\mid \abs{h(z)}\leq \sup \abs{h}(B)\text{ for any holomorphic }h \text{ on }\Omega\}.
\end{equation*}
 $\Omega$ is said to be holomorphically convex, whenever any of its compact subsets has relatively compact holomorphically convex hull in $\Omega$. 
\end{definition}

The Cartan--Thullen theorem reads as follows, see \cite{Cartan1932} for the original paper and \cite[Theorem 2.10.3 and Theorem 2.10.4, p. 77]{Klimek1991}, \cite[Theorem 3.4.5., p. 153]{Krantz2001} for more recent expositions.

\begin{theorem}\label{thm:cartan}
Suppose that $\Omega\subset\mathbb{C}^n$ is an open set. Then the following conditions are equivalent:
\begin{enumerate}[label=(\roman*), series=l_after]
\item\label{i:holoconv} for any compact set $K\subset\Omega$, its closed holomorphically convex hull is a compact subset of $\Omega$,
\item\label{i:domain} $\Omega$ is a domain of holomorphy, i.e., there are no open sets $\Omega_1,\Omega_2\subset\mathbb{C}^n$ such that 
\begin{enumerate}[label=(\alph*)]
\item\label{i:em} $\emptyset\neq \Omega_1\subset\Omega_2\cap \Omega$,
\item\label{i:in} $\Omega_2$ is connected and $\Omega_2\setminus\Omega\neq\emptyset$,
\item for each holomorphic $a$ on $\Omega$ there exists a holomorphic function $\hat{a}$ on $\Omega_2$ such that $\hat{a}=a$ on $\Omega_1$,
\end{enumerate}
\item\label{i:existence} $\Omega$ is a domain of existence, i.e., there exists holomorphic $a$ on $\Omega$ such that there are no open sets $\Omega_1,\Omega_2$ satisfying \ref{i:em}, \ref{i:in} and such that $a=\hat{a}$ on $\Omega_1$ for some holomorphic function $\hat{a}$ on $\Omega_2$.
\end{enumerate}
\end{theorem}

The Levi problem, posed by Levi in 1911 \cite{Levi1911}, was concerned with the following question: is an open set $\Omega\subset\mathbb{C}^n$ such that any of its boundary points admit a neighbourhood $U$ such that $U\cap\Omega$ is holomorphically convex, necessarily holomorphically convex itself?
The positive resolution of the problem is equivalent to demonstrating that the holomorphic convexity of $\Omega$ is equivalent to the following condition:
\vspace{0.3cm}

\begin{enumerate}[label=(\roman*), resume*=l_after]
    \item \label{i:ex}\emph{there exists a continuous, proper plurisubharmonic function} $p\colon\Omega\to\mathbb{R}$.
\end{enumerate}
\vspace{0.3cm}
The problem  remained unresolved for more than forty years until an affirmative resolution by Oka in 1950 \cite{Oka1950}, who used homological methods. Independent solutions were provided by Bremermann \cite{Bremermann1954} and by Norguet \cite{Norguet1954}.
Another approach, developed by H\"ormander \cite{Hormander1965}, see also \cite[Chapter 5]{Krantz2001}, employed methods of partial differential equations.
 We shall provide a simple proof of an analogue of the equivalence of \ref{i:holoconv}-\ref{i:ex} in the context of generalised convexity. In particular, we prove a generalisation of Theorem \ref{thm:cartan}.
 
The reader is referred to a survey on the Levi problem \cite{Siu1978}, where the Cartan--Thullen theorem is also presented \cite[Theorem 1.3., p. 482]{Siu1978}. For a modern discussion of the problem we refer to \cite[Theorem 3.3.5., p. 144, Theorem 3.4.5., p. 153]{Krantz2001}. 
A more recent paper on the topic is \cite{Sibony2018}.  

Using our abstract formulation of the Levi problem we show that the classical lemma of Bremerman and Lelong \cite[Theorem 2]{Bremermann1956} is equivalent to the affirmative answer to the Levi problem.

\subsection{Generalised convexity}

In \cite{Fan1963} Fan introduced a notion of generalised convexity with respect to a family of functions for closed sets. It allowed for a generalisation of  the Krein--Milman theorem.
Here we shall be concerned with sets that are convex with respect to a cone  $\mathcal{F}$ of continuous functions on a topological space $\Omega$. A closed set $K\subset \Omega$ shall be called $\mathcal{F}$-convex whenever it is equal to its closed $\mathcal{F}$-convex hull
\begin{equation*}
\mathrm{clConv}_{\mathcal{F}}K=\{\omega\in\Omega\mid f(\omega)\leq \sup f(K)\text{ for all }f\in\mathcal{F}\}.
\end{equation*}

Not only does the notion of generalised convexity  generalise classical convexity, when $\mathcal{F}$ is the cone of lower semi-continuous and convex functions, but it also includes the notion of holomorphic convexity, when $\mathcal{F}$ is chosen to be the space of real-parts of holomorphic functions. 

The reader is referred to the book of H\"ormander \cite{Hormander2007} for an account comprising matter on notions of convexity with respect to: subharmonic functions, \cite[Chapter III]{Hormander2007}, plurisubharmonic  functions, \cite[Chapter  IV]{Hormander2007}. The latter is equivalent to the notion of  holomorphic convexity of complex analysis, cf. \cite[Chapter 3]{Krantz2001}.

We shall provide a notion of convexity of a set with respect to a family of functions, regardless of the closedness of the set. This is achieved by means of embedding the space via the Gelfand transform, see Section \ref{s:gelfand} and Section \ref{s:convex}.

The notion of convexity has been developed and further generalised over the years, see \cite{Dolecki1978}, \cite{Rubinov2000}, and \cite{Singer1997}. These developments were mainly concerned with various dualities,  notions of subdifferentials and Fenchel transforms, which  are parallel to the current developments.

\subsubsection{Examples}

In the book \cite{Hormander2007} H\"ormander exhibits several examples that show that the equivalence of \ref{i:holoconv} of Theorem \ref{thm:cartan} and the Levi-type condition \ref{i:ex} holds true in a variety of settings that involve generalised convexity.
These equivalences are:
\begin{enumerate}[label=(\roman*)]
\item \cite[Corollary 2.1.26., p. 59]{Hormander2007} shows that for an open set $U\subset \mathbb{R}^n$ its convexity is equivalent to the existence of a convex, proper function $p\colon U\to\mathbb{R}$,
 \item \cite[Theorem 3.2.31., p. 165]{Hormander2007} shows that if $U\subset\mathbb{R}^n$ is an open set, then the following conditions are equivalent:
 \begin{enumerate}
 \item there is a proper, continuous subharmonic function $p\colon U\to\mathbb{R}$,
 \item for any compact set $K\subset  U$ its subharmonic convex hull 
 \begin{equation*}
     \{\omega\in U\mid a(\omega)\leq \sup a(K)\text{ for all subharmonic }a\text{ on }U\}
 \end{equation*}
 is a compact subset of $U$,
 \end{enumerate}
\item \cite[Theorem 4.1.19., p. 236]{Hormander2007} shows that if $U\subset\mathbb{C}^n$ is an open set, then the following conditions are equivalent:
 \begin{enumerate}
 \item there is a proper, continuous plurisubharmonic function $p\colon U\to\mathbb{R}$, 
 \item for any compact set $K\subset  U$ its plurisubharmonic convex hull  
 \begin{equation*}
     \{\omega\in U\mid a(\omega)\leq \sup a(K)\text{ for all plurisubharmonic }a\text{ on }U\}
 \end{equation*}
 is a compact subset of $U$.
\end{enumerate}
\end{enumerate}

We shall see that these results of H\"ormander are instances of general Theorem \ref{thm:holoalgin}, in which we consider cones $\mathcal{F}$ of functions that satisfy the maximum principle. Another related result is Theorem \ref{thm:holoin} concerned with complete lattices generated by a linear subspaces of functions, see Definition \ref{def:lattice}.

We mention also the work of Range \cite{Range2012} presenting relations between pseudoconvexity and convexity.

We shall also see that the Cartan--Thullen theorem,  Theorem \ref{thm:cartan}, admits a generalisation to the setting of generalised convexity, under the assumption that the cone $\mathcal{F}$ is generated by a linear space of functions, see Definition \ref{def:lattice}, that is, any function in the cone is the supremum of a collection of functions in the linear space.

\subsection{Domains of holomorphy and domains of existence}

Below  novel definitions are provided  serving  as analogues of the notions of domain of holomorphy and domain of existence. These are Definition \ref{def:comp}, \ref{i:cauchy}, and Definition \ref{def:aspace} respectively.

\begin{definition}
Let $\Omega$ be a set and let $\mathcal{A}$ be a set of functions on $\Omega$. The coarsest topology on $\Omega$ with respect to which all functions in $\mathcal{A}$ are continuous we shall call the topology generated by $\mathcal{A}$ and denote by $\tau(\mathcal{A})$.
\end{definition}

\begin{remark}
The topology $\tau(\mathcal{A})$ defined above is generated by the basis sets of the form
\begin{equation*}
\{\omega\in\Omega\mid \abs{a_i(\omega)-a_i(\omega_0)}\leq \epsilon_i\text{ for }i=1,2,\dotsc,k\},
\end{equation*}
where $\epsilon_i>0$, $\omega_0\in\Omega$ and $a_i\in\mathcal{A}$ for $i=1,2,\dotsc,k$.
\end{remark}

\begin{definition}\label{def:comp}
Let $\Omega$ be a topological space and let $\mathcal{A}$ be a linear space of continuous  functions on $\Omega$. We shall write that:
\begin{enumerate}[label=(\roman*)]
\item a net $(\omega_{\alpha})_{\alpha\in A}$ in $\Omega$ is a Cauchy net with respect to $\mathcal{A}$ whenever for each $a\in\mathcal{A}$, $(a(\omega_{\alpha}))_{\alpha\in A}$ is a Cauchy net in $\mathbb{R}$,
\item\label{i:cauchy} $\Omega$ is complete with respect to $\mathcal{A}$ whenever every Cauchy net with respect to $\mathcal{A}$ in $\Omega$ is convergent in $\tau(\mathcal{A})$ to an element in $\Omega$.
\end{enumerate}
\end{definition}

\begin{definition}\label{def:aspace}
 Let $\Omega$ be a topological space and let  $\mathcal{A}$ be a linear space of functions on $\Omega$. Let $d$ be a metric on $\Omega$ that generates the same topology as $\tau(\mathcal{A})$.

We shall say that $\Omega$ is an $(\mathcal{A},d)$-space whenever there exists a \emph{single function} $a\in\mathcal{A}$ for which 
\begin{equation*}
    \lim_{\alpha\in A} a(\omega_{\alpha})=\infty
\end{equation*}
for any Cauchy net $(\omega_{\alpha})_{\alpha\in A}$ with respect to $d$ in $\Omega$ that does not converge in $\tau(\mathcal{A})$ to an element in $\Omega$.
\end{definition}

\begin{remark}
    Note that it is not always true that a function $a\in\mathcal{A}$ will extend to a continuous function on the completion of $\Omega$. However, it is so, if, for example, we assume that $a$ is uniformly continuous.
\end{remark}

\begin{example}\label{exa:holo}
Consider an open set $\Omega\subset\mathbb{C}^n$  and take $\mathcal{A}$ to be the space of real-parts of holomorphic functions on $\Omega$. If $\Omega$ is complete with respect to $\mathcal{A}$, then it is a domain of holomorphy. Indeed, suppose there exist  open sets $\Omega_1,\Omega_2\subset\mathbb{C}^n$ such that $\emptyset\neq \Omega_1\subset \Omega_2\cap\Omega$, $\Omega_2$ is connected, $\Omega_2\setminus\Omega\neq\emptyset$,  with the property that for any holomorphic $a$ on $\Omega$, its restriction to $\Omega_1\cap\Omega$ extends to a holomorphic function on $\Omega_2$. Let $\omega\in\partial\Omega\cap \Omega_2$. By the extension properties of the holomorphic functions on $\Omega$, any sequence $(\omega_n)_{n=1}^{\infty}\subset\Omega_2$ converging to $\omega$, is a Cauchy sequence with respect to $\mathcal{A}$. However, it is not convergent in $\Omega$, so $\Omega$ is not complete with respect to $\mathcal{A}$. 

Similarly, if $\Omega$ is an $(\mathcal{A},d)$-space, where $d$ is the standard Euclidean metric, then it is a domain of existence. 
Let $a\in\mathcal{A}$ be a function such that
\begin{equation*}
    \lim_{\alpha\in A} a(\omega_{\alpha})=\infty
\end{equation*}
for any Cauchy net $(\omega_{\alpha})_{\alpha\in A}$ with respect to $d$ in $\Omega$ that does not converge in $\tau(\mathcal{A})$.
Then there are no open sets $\Omega_1,\Omega_2\subset\mathbb{C}^n$ such that $\emptyset\neq \Omega_1\subset \Omega_2\cap\Omega$, $\Omega_2$ is connected, $\Omega_2\setminus\Omega\neq\emptyset$,  with the property that the restriction of $a$ to $\Omega\cap \Omega_1$ extends to a holomorphic function on $\Omega_2$. 
\end{example}

\subsection{Main results}

\begin{definition}
We shall say that a family of functions is \emph{symmetric} if together with any function it contains its negative.
\end{definition}

\begin{definition}
Let $\mathcal{F}$ be a cone of functions on a set $\Omega$. We say that $P$ is an $\mathcal{F}$-polytope whenever there exists a finite family $(f_i)_{i=1}^k$ of functions in $\mathcal{F}$ such that
\begin{equation*}
P=\big\{\omega\in\Omega\mid f_i(\omega)\leq 1\text{ for  }i=1,2\dotsc,k\big\}.
\end{equation*}
We shall say that an $\mathcal{F}$-polytope is symmetric whenever the corresponding family can be taken to be symmetric.
\end{definition}

In the first of our main results, we provide an abstract version of the Levi problem, which is the equivalence of \ref{i:fconstantin}-\ref{i:exhaustionbound} below. 
Furthermore, under the additional assumption of separability of $\Omega$, we prove an abstract generalisation of the Cartan--Thullen theorem, cf. Theorem \ref{thm:cartan}: we show that conditions \ref{i:fconstantin}-\ref{i:exhaustionbound} are equivalent to $\Omega$ being complete with respect to $\mathrm{cl}\mathcal{A}$, which is an analogous notion to the notion of domain of holomorphy. If $\mathcal{A}$ consists of real-parts of a complex algebra, then these conditions are also equivalent to $\Omega$ being an $(\mathrm{cl}\mathcal{A},d)$-space, for any $d$ that metrises $\Omega$, which is an analogue of the notion of domain of existence. 

\begin{theorem}\label{thm:holoin}
Let $\Omega$ be a $\sigma$-compact, locally compact Hausdorff topological space.  Let $\mathcal{A}$ be a linear space of continuous functions on $\Omega$ that contains constants and separates points of $\Omega$.
The following conditions are equivalent:
\begin{enumerate}[label=(\roman*)]
\item\label{i:fconstantin} for any compact set $K\subset\Omega$ and for any $C\geq 1$ the set
\begin{equation*}
\{\omega\in\Omega\mid a(\omega)\leq C\sup \abs{a}(K)\text{ for all }a\in\mathcal{A}\}
\end{equation*}
is a compact subset of $\Omega$,
\item\label{i:exhaustionfin}
there exists a non-negative and proper function $p$ on  $\Omega$ such that  
\begin{equation*}
p=\sup\{a_{\alpha}\mid \alpha\in A\}
\end{equation*}
for some symmetric family $(a_{\alpha})_{\alpha\in A}$
of functions in $\mathcal{A}$, such that for any compact $K\subset\Omega$ there exists $\{\alpha_1,\dotsc,\alpha_k\}\subset A$ such that
\begin{equation*}
p=\sup\{\abs{a_{\alpha_i}}\mid i=1,\dotsc,k\}
\end{equation*}
on $K$,
\item\label{i:exhaustion} there exists a non-negative, continuous and proper function $p$ on  $\Omega$ such that  
\begin{equation*}
p=\sup\{a_{\alpha}\mid \alpha\in A\}
\end{equation*}
for some symmetric family $(a_{\alpha})_{\alpha\in A}$
of functions in $\mathcal{A}$,
\item\label{i:exhaustionbound} there exists a non-negative, proper function $p$ on  $\Omega$, bounded on compacta, such that  
\begin{equation*}
p=\sup\{a_{\alpha}\mid \alpha\in A\}
\end{equation*} for some symmetric family $(a_{\alpha})_{\alpha\in A}$ of functions in $\mathcal{A}$.
\end{enumerate}
  If $\mathcal{A}$  consists of real-parts of a complex algebra, then \ref{i:fconstantin} is equivalent to the existence of a family $(P_i)_{i=1}^{\infty}$ of compact, symmetric $\mathcal{A}$-polytopes such that
\begin{equation*}
\bigcup_{i=1}^{\infty}P_i=\Omega, P_i\subset\mathrm{int}P_{i+1}\text{  for }i=1,2,\dotsc.
\end{equation*} 

If $\Omega$ is separable then \ref{i:fconstantin} is equivalent to $\Omega$ being complete with respect to $\mathrm{cl}\mathcal{A}$, where the closure is taken with respect to the compact-open topology. 

 If $\Omega$ is separable and $\mathcal{A}$ consists of real-parts of a complex algebra and $d$ is any metric on $\Omega$ that generates the same topology as $\tau(\mathcal{A})$ , then \ref{i:fconstantin} is equivalent to $\Omega$ being a $(\mathrm{cl}\mathcal{A},d)$-space.
\end{theorem}

\begin{remark}
If $\Omega$ is $\sigma$-compact, locally compact, Hausdorff topological space and $\mathcal{A}$  consists of continuous functions that separate points of $\Omega$, then the topology of $\Omega$ is the weak topology $\tau(\mathcal{A})$ induced by $\mathcal{A}$, as we shall see in  Lemma \ref{lem:embed}.
\end{remark}

\begin{remark}
If $\Omega$ is $\sigma$-compact, locally compact, Hausdorff topological space and $\mathcal{A}$ is a closed algebra of real-valued continuous functions that separates points of $\Omega$, then $\mathcal{A}$ is the space of all continuous functions on $\Omega$, as follows by the Stone--Weierstrass theorem. Therefore we do not consider in the theorem the analogues of the  statements for real algebras that were stated for complex algebras.
\end{remark}

\begin{remark}
The function  $p$ of \ref{i:exhaustionfin}, \ref{i:exhaustion}  or of \ref{i:exhaustionbound} can be thought of as an analogue of a norm on the space $\Omega$.
If $\mathcal{A}$ is the space of continuous linear functions on a finite-dimensional linear space $\Omega$, then $p$ is a genuine norm on $\Omega$.
\end{remark}

\begin{remark}\label{rem:psh}
Let us compare condition \ref{i:fconstantin} of Theorem \ref{thm:holoin}  and the condition \ref{i:holoconv} of Theorem \ref{thm:cartan}. Observe that \ref{i:fconstantin} involves arbitrary constants $C\geq 1$. However, we shall in Proposition \ref{pro:fconvexity}, that when $\mathcal{A}$ consists of real-parts of a complex algebra $\mathcal{B}$, then 
\begin{equation*}
\{\omega\in \Omega\mid a(\omega)\leq C\sup \abs{a}(K)\text{ for  all }a\in\mathcal{A}\}=\mathrm{clConv}_{\mathcal{A}}K
\end{equation*}
for any compact $K\subset \Omega$ and any $C\geq 1$.

Let us,  however, observe that the modulus, which appears in \ref{i:fconstantin}, does play a r\^ole in the general case. Indeed, let us fix a compact set $K\subset\Omega$. Proposition \ref{pro:convexequivalent} shows that for $C\geq 1$ the sets
\begin{equation*}
\{\omega\in \Omega\mid a(\omega)\leq C\sup a(K)\text{ for  all }a\in\mathcal{A}\}
\end{equation*}
are all equal. In particular, they are equal to $\mathrm{clConv}_{\mathcal{A}}K$.
However, note that in general the sets $\mathrm{clConv}_{\mathcal{A}}K$ and $\{\omega\in \Omega\mid a(\omega)\leq C\sup \abs{a}(K)\text{ for  all }a\in\mathcal{A}\}$ are not equal, as exemplified by the space $\mathcal{A}$ of affine functions on an open, proper convex subset $\Omega$ a finite-dimensional space,   see Remark \ref{rem:nonequi}.
\end{remark}

\subsubsection{Cones satisfying the maximum principle}

The second main result is concerned with another abstract version of the Levi problem in the generalised convexity setting. Instead of considering linear spaces of functions, we are now investigating the case of convexity with respect to a convex cone of functions that satisfies the maximum principle.

\begin{definition}
We shall say that a convex cone $\mathcal{F}$ of continuous functions on a topological space $\Omega$ satisfies the maximum principle whenever for any:
\begin{enumerate}[label=(\roman*)]
\item function $f\in\mathcal{F}$ that is non-constant,
\item  compact set $K\subset\Omega$,
\item open set $U\subset \Omega$ such that $ K\subset  U$, 
\end{enumerate} 
there is
\begin{equation*}
\sup f(K)<\sup f(U).
\end{equation*}
\end{definition}

\begin{remark}
If  $\Omega$ is compact, then one may take $K=\Omega$. Thus, if $\Omega$ is compact and if $\mathcal{F}$ satisfies the maximum principle, then any function $f\in\mathcal{F}$ is constant.
\end{remark}

We will show that the existence of a continuous $\mathcal{F}$-exhaustion function $p$, cf. Definition \ref{def:exhaustion}, is equivalent to compactness of $\mathrm{clConv}_{\mathcal{F}}K$ for $K\subset\Omega$ compact, under the assumption that the cone  $\mathcal{F}$ of considered functions satisfies the maximum principle.  
\begin{theorem}\label{thm:holoalgin}
Let $\Omega$ be a $\sigma$-compact, locally compact, Hausdorff connected topological space.  Let $\mathcal{F}$ be a convex cone of continuous functions on $\Omega$ that contains constants and satisfies the maximum principle.  
Then the following  conditions are equivalent:
\begin{enumerate}[label=(\roman*)]
\item \label{i:foneRin} for any compact set $K\subset\Omega$ the set
\begin{equation*}
\mathrm{clConv}_{\mathcal{F}}K=\{\omega\in\Omega\mid f(\omega)\leq \sup f(K)\text{ for all }f\in\mathcal{F}\}
\end{equation*}
is a compact subset of $\Omega$,
\item\label{i:exhaustionfincon}
there exists a non-negative and proper function $p$ on  $\Omega$ such that  
\begin{equation*}
p=\sup\{f_{\alpha}\mid \alpha\in A\}
\end{equation*}
for some family $(f_{\alpha})_{\alpha\in A}$
of functions in $\mathcal{F}$, such that for any compact $K\subset\Omega$ there exists $\{\alpha_1,\dotsc,\alpha_k\}\subset A$ such that
\begin{equation*}
p=\sup\{f_{\alpha_i}\mid i=1,\dotsc,k\}
\end{equation*}
on $K$,
\item\label{i:exhaustionconR} there exists a non-negative, continuous and proper function $p$ on $\Omega$
such that  
\begin{equation*}
p=\sup\{f_{\alpha}\mid \alpha\in A\}
\end{equation*}
for some family $(f_{\alpha})_{\alpha\in A}$
of functions in $\mathcal{F}$,
\item\label{i:exhaustionboundR}
there exists a non-negative, proper function $p$ on  $\Omega$, bounded on compacta, such that  
\begin{equation*}
p=\sup\{f_{\alpha}\mid \alpha\in A\}
\end{equation*} for some family $(f_{\alpha})_{\alpha\in A}$ of functions in $\mathcal{F}$,
\item\label{i:compactaf}
there exists a family $(P_i)_{i=1}^{\infty}$ of compact $\mathcal{F}$-polytopes such that
\begin{equation*}
\bigcup_{i=1}^{\infty}P_i=\Omega, F_i\subset\mathrm{int}P_{i+1}\text{  for }i=1,2,\dotsc.
\end{equation*} 
\end{enumerate}
\end{theorem}

\begin{remark}
Suppose that a linear space $\mathcal{A}$ of continuous functions on $\Omega$ satisfies the maximum principle. Then we may apply both Theorem \ref{thm:holoin} and Theorem \ref{thm:holoalgin}, with $\mathcal{F}=\mathcal{A}$. We see that the conditions of Theorem \ref{thm:holoalgin} are implied by  the conditions of Theorem \ref{thm:holoin}.
\end{remark}

\begin{remark}\label{rem:alg}
As we have already  noted in Remark \ref{rem:psh}, when $\mathcal{F}$ consists of real-parts of a complex algebra, then the conditions \ref{i:fconstantin} of Theorem \ref{thm:holoin}, where $\mathcal{A}=\mathcal{F}$, and \ref{i:foneRin}  of Theorem \ref{thm:holoalgin} are equivalent, see Proposition \ref{pro:fconvexity}. For a characterisation of \ref{i:foneRin}  of Theorem \ref{thm:holoalgin}, see Proposition \ref{pro:convexequivalent}. 
\end{remark}

\begin{definition}
We  shall say that the cone $\mathcal{F}$ of functions on a topological space $\Omega$ is local, whenever for any family $(U_i)_{i\in I}$ of distinct connected components of $\Omega$ and any family $(f_i)_{i\in I}\subset\mathcal{F}$, there exists $f\in\mathcal{F}$ such that
\begin{equation*}
f=f_i\text{ on }U_i.
\end{equation*}
\end{definition}

\begin{remark}
If the cone $\mathcal{F}$ is local, then we may prove Theorem \ref{thm:holoalgin} dropping the assumption on connectedness of $\Omega$. Indeed, $\sigma$-compactness of $\Omega$ implies that it has at most countable many connected components $(U_i)_{i=1}^{\infty}$. We apply Theorem \ref{thm:holoalgin} on each of the components separately. For each of the components we pick $p_i\in\mathcal{F}$ that satisfies \ref{i:exhaustionboundR} on $U_i$. Then we set 
\begin{equation*}
p=p_i+i\text{ on }U_i \text{ for all }i=1,2,\dotsc.
\end{equation*}
Since $\mathcal{F}$ is local, we see that $p$ satisfies \ref{i:exhaustionboundR} for the space $\Omega$. The other implications of the generalisation of Theorem \ref{thm:holoalgin} follow as in the proof of Theorem \ref{thm:holoalgin}. 
\end{remark}

\subsubsection{Equivalence of the Bremermann--Lelong lemma and the positive resolution of the Levi problem}

In the field of complex analysis of several complex variables, the lemma of Bremermann and of Lelong \cite[Theorem 2]{Bremermann1956} states that if $\Omega\subset\mathbb{C}^n$ is a domain of holomorphy, then any plurisubharmonic function on $\Omega$ is a Hartogs function. The set of Hartogs functions on $\Omega$, cf. \cite[Section 2]{Bremermann1956}, is the smallest convex cone of functions on $\Omega$ that contains the logarithms of moduli of holomorphic functions on $\Omega$ and is closed under taking:
\begin{enumerate}[label=(\roman*)]
    \item locally bounded suprema,
    \item monotonically decreasing limits.
\end{enumerate}
Moreover, the function belongs to the set of Hartogs functions on $\Omega$ if and only if its restrictions to any open, precompact subset $U$ of $\Omega$ is a Hartogs function on $U$, and any upper semi-continuous regularisation of a Hartogs function is again a Hartogs function. 

We shall show in Theorem \ref{thm:bremermann} that our methods allow to prove that the positive resolution of the Levi problem is equivalent to the lemma of Bremermann and Lelong.

\subsection{Examples}

We shall now exhibit several examples of how Theorem \ref{thm:holoin} and Theorem \ref{thm:holoalgin} apply to various settings. 

\subsubsection{Holomorphic convexity}

The example that inspired this research is concerned with holomorphic convexity.

The following theorem, which implies Theorem \ref{thm:cartan}, follows from Theorem \ref{thm:holoin}.

\begin{theorem}\label{thm:levi}
Suppose that $\Omega\subset\mathbb{C}^n$ is an open set and let $\mathcal{B}$ denote the space of holomorphic functions on $\Omega$. Then the following are equivalent:
\begin{enumerate}[label=(\roman*)]
\item\label{i:compa} any compact set $K\subset\Omega$, the set
\begin{equation*}
\{\omega\in\Omega\mid \abs{b(\omega)}\leq \sup \abs{b}(K)\text{ for any }b\in\mathcal{B}\}
\end{equation*}
 is a compact subset of $\Omega$,
  \item\label{i:bcopen} $\Omega$ is complete with respect to $\mathcal{B}$.
   \item\label{i:exho} there exists a non-negative, continuous and proper function $p$ on $\Omega$ such that
 \begin{equation*}
 p=\sup\{\mathfrak{Re}b_{\alpha}\mid \alpha\in A\}
 \end{equation*}
 for some family $(b_{\alpha})_{\alpha\in A}\subset\mathcal{B}$,
 \item\label{i:polytopes} there exists a family  $(P_i)_{i=1}^{\infty}$ of holomorphic polytopes such that
 \begin{equation*}
 \bigcup_{i=1}^{\infty}P_i=\Omega\text{ and  }P_i\subset\mathrm{int}P_{i+1}\text{ for }i=1,2,\dotsc,
 \end{equation*}
 \item\label{i:bdopen} $\Omega$ is a $(\mathcal{B},d)$-space, where $d$ is any metric on $\Omega$ that generates the standard topology.
  \end{enumerate}
\end{theorem}

Let us now discuss the relationship between the  above theorem and Theorem \ref{thm:holoin}.

The condition that $\Omega$ is a $(\mathcal{B},d)$-space readily implies that it is a domain of existence, see Example \ref{exa:holo}, which in turn trivially implies that it is a domain of holomorphy.
Similarly, the condition that $\Omega$ is complete with respect to $\mathcal{B}$ shows that it is a domain of holomorphy. This is known to imply \ref{i:holoconv} of Theorem \ref{thm:cartan} by a simple argument involving  power series expansion, see \cite[Theorem 3.4.5, p. 153]{Krantz2001}. Let us stress that this argument is not needed to prove the equivalence in Theorem \ref{thm:levi} above -- it is only needed if we want to prove the original version of the Cartan--Thullen theorem.
Moreover, the implications in our proof rely on the theorems of functional analysis: the Banach--Alaoglu theorem and the Banach--Steinhaus uniform boundedness theorem. Our  proof thus allows to locate the classical Cartan--Thullen theorem in a broader context -- as an instance of a general result in generalised convexity.

 Moreover, Theorem \ref{thm:levi} also implies the modification of \ref{i:ex} to \ref{i:exho} of Theorem \ref{thm:levi}. 

Let us note, that we may as well employ Theorem \ref{thm:holoalgin} directly to the setting of an open subset of $\mathbb{C}^n$ and the cone of plurisubharmonic functions. Then we recover also
\cite[Theorem 4.1.19., p. 236]{Hormander2007}.

\subsubsection{Stein spaces}

Let us recall the following definition.

\begin{definition}\label{def:stein}
A complex analytic space $\Omega$ is said to be a Stein space whenever:
\begin{enumerate}[label=(\roman*)]
\item the space $\mathcal{B}$ of holomorphic functions on $\Omega$ separates points of $\Omega$,
\item for any compact set $K\subset\Omega$, the set
\begin{equation*}
\{\omega\in\Omega\mid \abs{b(\omega)}\leq \sup \abs{b}(K)\text{ for any }b\in\mathcal{B}\}
\end{equation*}
 is a compact subset of $\Omega$.
\end{enumerate}
\end{definition}

Theorem \ref{thm:holoin}, together with Proposition \ref{pro:fconvexity} implies the following characterisation of Stein spaces.

\begin{theorem}\label{thm:stein}
Suppose that $\Omega$ is $\sigma$-compact, locally compact, separable complex analytic space such that the space $\mathcal{B}$ of holomorphic functions on $\Omega$ separates points of $\Omega$. Then the following are equivalent:
\begin{enumerate}[label=(\roman*)]
\item $\Omega$ is a Stein space,
\item any compact set $K\subset\Omega$, the set
\begin{equation*}
\{\omega\in\Omega\mid \abs{b(\omega)}\leq \sup \abs{b}(K)\text{ for any }b\in\mathcal{B}\}
\end{equation*}
 is a compact subset of $\Omega$,
  \item\label{i:bc} $\Omega$ is complete with respect to $\mathcal{B}$.
   \item\label{i:exh} there exists a non-negative, continuous and proper function $p$ on $\Omega$ such that
 \begin{equation*}
 p=\sup\{\mathfrak{Re}b_{\alpha}\mid \alpha\in A\}
 \end{equation*}
 for some family of holomorphic functions $(b_{\alpha})_{\alpha\in A}\subset\mathcal{B}$,
  \item\label{i:polytopesstein} there exists a family  $(P_i)_{i=1}^{\infty}$ of holomorphic polytopes such that
 \begin{equation*}
 \bigcup_{i=1}^{\infty}P_i=\Omega\text{ and  }P_i\subset\mathrm{int}P_{i+1}\text{ for }i=1,2,\dotsc,
 \end{equation*}
 \item\label{i:bd} $\Omega$ is a $(\mathcal{B},d)$-space for any  metric $d$ that generates the same topology as $\tau(\mathcal{B})$.
  \end{enumerate}
\end{theorem}

Let us recall that, similarly as in the case of open subsets of $\mathbb{C}^n$, there exists a characterisation of Stein manifolds that relies on the existence of strictly plurisubharmonic exhaustion function, see \cite[Theorem 5.1.6., p. 117]{Hormander1990}.

Let us stress that, to the best of our knowledge, the conditions \ref{i:bc}, \ref{i:exh} and \ref{i:bd} of Theorem \ref{thm:stein} are new. However, we remark that \cite[Theorem 5.4.2., p. 138]{Hormander1990}  is similar to the completeness condition and \cite[Lemma 5.3.7., p. 133]{Hormander1990} implies \ref{i:polytopesstein}. The condition \ref{i:bc} shows that the intention of Stein, the founder of the notion, to name  spaces that satisfy the conditions of Definition \ref{def:stein} as holomorphically complete, see \cite{Stein1951}, was accurate.

\subsubsection{Plurisubharmonic functions}

When $\mathcal{F}$ is the cone of plurisubharmonic functions, then Theorem \ref{thm:holoalgin} reproves \cite[Theorem 4.1.19., p. 236]{Hormander2007}.

\subsubsection{Subharmonic functions}

When $\mathcal{F}$ is the cone of subharmonic functions, then Theorem \ref{thm:holoalgin} reproves \cite[Theorem 3.2.31., p. 165]{Hormander2007}. 

\subsubsection{Convex functions}

 An analogous characterisation for the space of convex functions, see \cite[Theorem 2.1.25., p. 58]{Hormander2007} and its corollary \cite[ Corollary 2.1.26., p. 59]{Hormander2007}, can again be inferred from Theorem \ref{thm:holoalgin}.

\subsection{Gelfand transform and relation to convexity}

The notion of convexity that we study here is naturally linked to the classical convexity via the Gelfand transform $\Phi\colon\Omega\to\mathcal{A}^*$, which is defined by the  formula
\begin{equation*}
    \Phi(\omega)(a)=a(\omega)\text{ for all }\omega\in\Omega\text{ and }a\in\mathcal{A}.
\end{equation*}
In our setting it can be viewed as a way of linearisation of the considered problems. The linearisation technique has been successful in a great many contexts, e.g., in the context of Lipschitz-free spaces, see \cite{Weaver1999}.

Using the Gelfand transform, we can transfer linear notions to the setting of the space $\Omega$, in a similar way to transference of the notion of convex hull, see Section \ref{s:convex}, Definition \ref{def:conv}. This allows for a reduction of investigations of an abstract space $\Omega$ equipped with a linear space $\mathcal{A}$ of functions on $\Omega$ to studying subsets of a linear space and linear, continuous functionals on that subset.

\subsection{Methods of proofs -- functional analytic ingredient}

The main novelty in our approach lies in the use of the classical theorems of functional analysis. We equip space the $\mathcal{A}$ with the compact-open topology. This is a topology generated by the family of semi-norms
\begin{equation*}
    \norm{a}_{\mathcal{C}(K)}=\sup\{\abs{a(\omega)}\mid\omega\in K\},
\end{equation*}
where $K$ are compact subsets of $\Omega$. Let us assume that $\mathcal{A}$ is closed. The sets
\begin{equation}\label{eqn:formm}
    \{\omega\in\Omega\mid a(\omega)\leq C\sup\abs{a}(K)\},
\end{equation}
considered in \ref{i:fconstantin} of Theorem \ref{thm:holoin}  for compact $K\subset\Omega$ and $C\geq 0$, consists precisely of these continuous functionals on $\mathcal{A}$ that are continuous with respect to the semi-norm $\norm{\cdot}_{\mathcal{C}(K)}$, and which are evaluations at the points of $\Omega$. This is to say 
\begin{equation*}
 \{\omega\in\Omega\mid \{a(\omega)\leq C\sup\abs{a}(K)\}=   \Phi^{-1}(R\cap\Phi(\Omega)),
\end{equation*}
where $R$ is a polar set of a neighbourhood of zero in $\mathcal{A}$. The Banach--Alaoglu theorem shows that $R$ is a compact set in $\mathcal{A}^*$. Note that $\Phi$ is a homeomorphism onto its image. Moreover, completeness of $\Omega$ with respect to $\mathcal{A}$ is equivalent to the completeness of $\Phi(\Omega)\subset\mathcal{A}^*$, and thus also  equivalent to its closedness. Thus, we see that completeness implies that the sets of the form (\ref{eqn:formm}) are compact.

The converse implication is a consequence of the Banach--Steinhaus uniform boundedness principle. Indeed,   a Cauchy sequence of elements in $\Omega$ induces the sequence of functionals in $\mathcal{A}^*$ -- evaluations at the points of the sequence. These functionals are pointwise bounded and therefore there is a neighbourhood of zero of $\mathcal{A}$ on which the sequence is bounded.  This means that the original Cauchy sequence belonged to a set of the form (\ref{eqn:formm}), whose compactness implies convergence of the considered sequence.

In Lemma \ref{lem:nonmetri} we present an alternative proof of one of the above implications, under an additional assumption that $\mathcal{A}$ consists of real-parts of a complex algebra. The technique we use in this alternative proof is not new -- however, it has not been applied previously in the context of generalised convexity.

The equivalence of \ref{i:fconstantin} of Theorem \ref{thm:holoin} and that $\Omega$ is an $(\mathrm{cl}\mathcal{A},d)$-space is an application of the technique usually employed in the setting of complex analysis of several variables. We include the proof for the sake of completeness. 

Another novelty in our developments is the proof of the equivalences of \ref{i:fconstantin}-\ref{i:exhaustionbound} in Theorem \ref{thm:holoin}. As we have mentioned already the previous approaches to the Levi problem used involved analytical or homological tools. In our approach, the desired function is constructed directly in an elementary way that relies on the $\sigma$-compactness assumption on the space $\Omega$. The construction in presented in detail in Lemma \ref{lem:exhaustionC}.

The proof of Theorem \ref{thm:holoalgin} follows similar steps to the proof of the equivalences of \ref{i:fconstantin}-\ref{i:exhaustionbound} in Theorem \ref{thm:holoin}. 

\subsubsection{Relation to theorems of functional analysis}
Let us note further analogies of Theorem \ref{thm:holoin} and Theorem  \ref{thm:holoalgin} with classical results of functional analysis, the statements of which can be found in \cite{Rudin1991}.
The equivalence of completeness and Theorem \ref{thm:holoin}, \ref{i:fconstantin} is an analogue of the Krein--\v{S}mulain theorem, Theorem \ref{thm:krein-smulian}.  Moreover, that  completeness of $\Omega$ implies  \ref{i:foneRin} of Theorem \ref{thm:holoalgin} follows from the Mazur theorem.

\subsection{Applications}

An application of the results obtained in this paper is concerned with a generalisation of martingale transport of measures, which we investigate in \cite{Ciosmak2024} and in \cite{Ciosmak20241}. Theorem \ref{thm:holoin} is used there to  provide a characterisation of spaces that admit an exhaustion function, whose existence is an assumption in a generalisation of the Strassen theorem proven in \cite{Ciosmak2024}. In the setting of such spaces in \cite{Ciosmak2024}  we extend the results of De March and Touzi \cite{Touzi2019}, Ob\l\'oj and Siorpaes \cite{Obloj2017}, Ghoussoub, Kim, Lim \cite{Ghoussoub2019} to the setting of martingales with respect to a linear space of functions. Let us also mention  \cite{Ciosmak2023}, where we study a related concept of ordering of measures with respect to a cone of functions.

The  aforementioned extension allows for a new type of localisation-type results. Another type of  localisation was studied in \cite{Ciosmak2021}. This approach is related to optimal transport of vector measures \cite{Ciosmak20211}, \cite{Ciosmak20212} and extensions of Lipschitz maps  \cite{Ciosmak20213, Ciosmak20242}.

The notions of generalised  convexity and of submartingale transport of measures, explored in \cite{Ciosmak20241}, have found its applications in a  study of  the monopolist's problem  by Figalli, Kim  and McCann \cite{Figalli2011}  and by McCann and Zhang in \cite{McCann2019}.

\subsection{Structure of the paper}
In Section \ref{s:prelim} we provide several definitions and lemmata.
In Section \ref{s:gelfand} we define the Gelfand transform and study its basic properties.
In Section \ref{s:convex} we provide definition of convexity with respect to a convex cone of functions.   Proposition \ref{pro:fconvexity} shows a characterisation of closed convex hulls with respect to a linear space of real-parts of a complex algebra.
In Section  \ref{s:proofs} we prove the main results, Theorem \ref{thm:holoin}, Theorem \ref{thm:holoalgin} and Theorem \ref{thm:levi}. We also elaborate on the analogy of Theorem \ref{thm:holoin} and the Krein--\v{S}mulian theorem.
In Section \ref{s:bremermann} we show that the lemma of Bremermann and Lelong implies the resolution of  the Levi problem.

\section{Preliminaries}\label{s:prelim}

Let us recall several definitions.

\begin{definition}
We  shall say  that a topological space $\Omega$ is \emph{exhaustible by compacta} if there exists a sequence $(K_i)_{i=1}^{\infty}$ of compact subsets of $\Omega$ such that
\begin{enumerate}[label=(\roman*)]
\item $\Omega=\bigcup_{i=1}^{\infty}K_i$,
\item $K_i\subset \mathrm{int}K_{i+1}$ for $i=1,2,\dotsc$.
\end{enumerate}
\end{definition}

\begin{definition}
We shall say that $\Omega$ is $\sigma$-compact whenever there exist compact subsets $(K_i)_{i=1}^{\infty}$ whose union is $\Omega$.
\end{definition}

\begin{remark}\label{rem:sigma}
It is readily seen that the condition that $\Omega$ is exhaustible by compacta is equivalent to  $\Omega$ being $\sigma$-compact and locally compact.\
\end{remark}

\begin{definition}\label{def:tauk}
Let $\mathcal{K}$ denote the family of all compact subsets of $\Omega$. For any compact set $K\in\mathcal{K}$ and $a\in\mathcal{A}$ we put
\begin{equation*}
\norm{a}_{\mathcal{C}(K)}=\sup\{\abs{a(\omega)}\mid \omega\in K\}.
\end{equation*}
The family of semi-norms $(\norm{\cdot}_{\mathcal{C}(K)})_{K\in\mathcal{K}}$ on $\mathcal{A}$ defines a locally convex topology on $\mathcal{A}$, which we shall denote by $\tau_{\mathcal{K}}$.
\end{definition}

\begin{remark}\label{rem:tauk}
Let us remark that the topology $\tau_{\mathcal{K}}$ coincides with the compact-open topology on $\mathcal{A}$. 

We shall denote by $\mathcal{A}^*$ the space of all linear functionals on $\mathcal{A}$ that are continuous with respect to the topology generated by the above family of semi-norms. 
We shall equip $\mathcal{A}^*$ with the weak* topology generated by $\mathcal{A}$.
\end{remark}

\begin{remark}\label{rem:separable}
If $\Omega$ is $\sigma$-compact, then $\mathcal{A}$ in $\tau_{\mathcal{K}}$ is metrisable, as the topology is induced by a countable separating family of semi-norms \cite[Remark 1.38, c), p. 28]{Rudin1991}.
If  moreover $\Omega$ is separable, then for each compact set $K\subset\Omega$ the space $\mathcal{C}(K)$ is separable, as follows by the Stone--Weierstrass theorem, so that $\mathcal{A}$ is separable. Therefore separability and $\sigma$-compactness of $\Omega$ imply that $\mathcal{A}^*$ in the weak* topology is metrisable on its compact subsets. 
\end{remark}

\begin{definition}\label{def:lattice}
Let $\Omega$ be a set. Let $\mathcal{A}$ be a set of functions on $\Omega$ with values in $(-\infty,\infty]$. A set $\mathcal{L}$ of functions on $\Omega$ is said to be:
\begin{enumerate}[label=(\roman*)]
\item a convex cone whenever $\alpha f+\beta g \in\mathcal{L}$ for any $f,g\in\mathcal{L}$ and any numbers $\alpha,\beta\geq 0$;
\item stable under suprema provided that the supremum $\sup\{f_i\mid i \in I\}$  of any family $(f_i)_{i\in I}\subset\mathcal{L}$ belongs to $\mathcal{L}$;
\item a complete lattice cone whenever it is a convex cone stable under suprema;
\item the complete lattice cone generated by $\mathcal{A}$ whenever it is the smallest complete lattice cone of functions on $\Omega$ containing $\mathcal{A}$;
\item separating points of $\Omega$ whenever for any two distinct $\omega_1,\omega_2\in\Omega$, there exists $f\in\mathcal{L}$ such that $f(\omega_1)\neq f(\omega_2)$.
\end{enumerate}
\end{definition}

Let us stress that we allow the functions to take the value $+\infty$. Note however that this is not allowed for functions in linear subspaces.

We refer to \cite{Fuchssteiner1981} for a study of the above notions in contexts related to this work.

\begin{lemma}\label{lem:supgen}
Suppose that $\mathcal{A}$ is a linear space of functions on $\Omega$. Let $\mathcal{F}$ denote the complete lattice cone generated by $\mathcal{A}$. Then for any $f\in\mathcal{F}$ there exists a family $(a_i)_{i\in I}\subset\mathcal{A}$ such that
\begin{equation}\label{eqn:forma}
f(\omega)=\sup\{a_i(\omega)\mid i\in I\}\text{ for all }\omega\in \Omega.
\end{equation}
\end{lemma}
\begin{proof}
Let $\mathcal{G}$ be the set of functions of the form (\ref{eqn:forma}).  We shall show that $\mathcal{G}$ is a convex cone. Let $\alpha,\beta\geq 0$, $f,g,\in\mathcal{G}$. Let $(a_i)_{i\in I}, (a_j)_{j\in J}\subset\mathcal{A}$ be families of functions in $\mathcal{A}$ corresponding to $f$ and $g$ respectively.
Then
\begin{equation*}
\alpha f+\beta g=\sup\{\alpha a_i+\beta a_j\mid j\in J, i\in I\}\in\mathcal{G}.
\end{equation*}
Clearly, $\mathcal{F}\supset\mathcal{G}$. Since $\mathcal{G}$ is a complete lattice cone, $\mathcal{F}=\mathcal{G}$.
\end{proof}

\begin{remark}\label{rem:convex}
If $Y$ is the space of all linear, continuous functionals on a linear topological space $X$, then the complete lattice cone generated by $Y$ is the set of all convex, lower semi-continuous functions on $X$. We refer the reader to \cite[Chapter 3, p. 13]{Phelps2001} for a proof of this fact.
\end{remark}

\begin{lemma}\label{lem:poussin}
Let $\mathcal{F}$ be the complete lattice cone of functions on $\Omega$ that contains  constants. 
Let $p\in\mathcal{F}$. Let $\xi\colon \mathbb{R}\to \mathbb{R}$ be convex and non-decreasing. Then the composition $\xi(p)$ belongs to $\mathcal{F}$.
\end{lemma}
\begin{proof}
As $\xi$ is convex and non-decreasing, there exists a family of non-decreasing affine functions $(\lambda_j)_{j\in J}$ on $\mathbb{R}$ such that for all $t\geq 0$ there is
\begin{equation*}
\xi(t)=\sup\{\lambda_j(t)\mid j\in J\}.
\end{equation*}
Then 
\begin{equation*}
\xi(p)(\omega)=\sup\{\lambda_j(p)(\omega)\mid j\in J\}.
\end{equation*}
For any $j\in J$, there exist $\alpha_j\geq 0$ and $\beta_j\in\mathbb{R}$ such  that $\lambda_j(t)=\alpha_j t+\beta_j$ for all $t\in\mathbb{R}$, so that $\alpha_j p+\beta_j\in\mathcal{F}$ for each $j\in  J$.
Thus $\xi(p)\in\mathcal{F}$.
\end{proof}

\begin{definition}\label{def:exhaustion}
Let $\Omega$ be a topological space. A map $p\colon \Omega\to \mathbb{R}$ is said to be proper if the preimage $p^{-1}(Z)$ ofr every compact sets in $Z\subset\mathbb{R}$ is  compact in $\Omega$.
Let $\mathcal{F}$ be a cone of functions on $\Omega$. We say that 
\begin{equation*}
p\colon \Omega\to [0,\infty)
\end{equation*}
is an $\mathcal{F}$-exhaustion function whenever it belongs to $\mathcal{F}$ and is proper.
\end{definition}

\section{Gelfand transform}\label{s:gelfand}

Given a topological space $\Omega$ and a linear space $\mathcal{A}$ of continuous functions on $\Omega$, one may consider the space of all linear functionals on $\mathcal{A}$ that are evaluations at points of $\Omega$.
Let $\Phi\colon\Omega\to\mathcal{A}^*$ be given by the formula $\Phi(\omega)(a)=a(\omega)$ for $\omega\in \Omega$ and $a\in\mathcal{A}$. We shall call $\Phi$ the Gelfand transform, in analogy with the notion studied in  the theory of Banach algebras, cf. \cite{Zelazko1973}.

\begin{lemma}\label{lem:embed}
Let $\Omega$ be a $\sigma$-compact and locally compact Hausdorff space and let $\mathcal{A}$ be a linear space of continuous functions on $\Omega$ that separates points of $\Omega$. Let $\mathcal{F}$ be the complete lattice cone generated by  $\mathcal{A}$. We equip $\mathcal{A}$ with the topology $\tau_{\mathcal{K}}$ and $\mathcal{A}^*$ with the weak* topology induced by $\mathcal{A}$. 
Then:
\begin{enumerate}[label=(\roman*)]
\item\label{i:homeo} $\Phi$ is a homeomorphism between $\Omega$ and $\Phi(\Omega)$,
\item\label{i:aextend} for any $a\in\mathcal{A}$, $a(\Phi^{-1})$ extends to a weakly* continuous linear map on $\mathcal{A}^*$, 
\item\label{i:arestrict} for any weakly* continuous linear map $h$ on $\mathcal{A}^*$, the function 
\begin{equation*}
\omega\mapsto h(\Phi(\omega))
\end{equation*}
belongs to $\mathcal{A}$,
\item\label{i:fextend} for any $f\in\mathcal{F}$, $f(\Phi^{-1})$ extends to a map in the complete lattice cone $G$  of functions on $\mathcal{A}^*$ generated by evaluations on elements of $\mathcal{A}$,
\item\label{i:frestrict} for any map $g\in G$, the function 
\begin{equation*}
\omega\mapsto g(\Phi(\omega))
\end{equation*}
belongs to $\mathcal{F}$.
\end{enumerate}
\end{lemma}
\begin{proof}

We shall first show that $\Phi$ is homeomorphism onto $\Phi(\Omega)$. As $\mathcal{A}$ separates points of $\Omega$ and consists of continuous functions, $\Phi$ is injective and continuous. Therefore, for any  compact $K\subset \Omega$, the restriction of $\Phi$ to $K$ is a homeomorphism onto the corresponding image.
As $\Omega$ is exhaustible by compacta, cf. Remark \ref{rem:sigma}, there exist compact sets $(K_i)_{i=1}^{\infty}$ such that
\begin{equation*}
\bigcup_{i=1}^{\infty}K_i=\Omega, K_i\subset\mathrm{int}K_{i+1}\text{ for }i=1,2,\dotsc.
\end{equation*}
Therefore for any open set $U\subset\Omega$ the set
\begin{equation*}
\Phi(U)=\bigcup_{i=1}^{\infty}\Phi(U\cap\mathrm{int}K_i)
\end{equation*}
is open. This is to say, $\Phi$ is a homeomorphism.

Let $a\in\mathcal{A}$. Extension in \ref{i:aextend} is given by the formula
\begin{equation*}
\mathcal{A}^*\ni a^*\mapsto a^*(a)\in\mathbb{R}.
\end{equation*}
We shall now define an extension for $f\in\mathcal{F}$. By Lemma \ref{lem:supgen} for $f\in\mathcal{F}$ there exists a family $(a_i)_{i\in I}$ of elements of $\mathcal{A}$ such that
\begin{equation}\label{eqn:repefi}
f=\sup\{a_i\mid i\in I\}.
\end{equation}
We define a map
\begin{equation*}
\mathcal{A}^*\ni a^*\mapsto \sup\{a^*(a_i)\mid i\in I\}\in (-\infty,\infty].
\end{equation*}
By (\ref{eqn:repefi}), this is an extension. As functionals $a^*\mapsto a^*(a)$ are weakly* continuous, this extension belongs to the complete lattice cone $G$ generated by $\mathcal{A}$. Thus \ref{i:fextend} is proven.

Note that any weakly* continuous linear functional $h$ on $\mathcal{A}^*$ is given by some element $a_0\in \mathcal{A}$. Therefore for $\omega\in\Omega$, 
\begin{equation*}
 h(\Phi(\omega))=a_0(\omega).
\end{equation*}
This proves \ref{i:arestrict}. Point \ref{i:frestrict} is proven similarly.
\end{proof}

\begin{remark}\label{rem:metrisable}
If $\Omega$ is exhaustible by compacta and separable,  $\mathcal{A}$ is closed in $\tau_{\mathcal{K}}$ and separates points of $\Omega$, then $\Omega$ is metrisable. By the Urysohn theorem it suffices to show that $\Omega$ has a countable basis of neighbourhoods and that it is completely regular. 
 
As $\Omega$ is a countable union of open sets which are metrisable, hence it admits  a countable basis of neighbourhoods. 
Lemma \ref{lem:embed} shows that the topology on $\Omega$ is the weak topology induced by $\mathcal{A}$. Thus, it is completely regular. 
\end{remark}

\section{$\mathcal{F}$-convex sets}\label{s:convex}

\begin{definition}\label{def:convexhull}
Let $\mathcal{F}$ be a convex cone of functions on $\Omega$. Let $S\subset\Omega$. We define closed $\mathcal{F}$-convex hull of $S$ by the formula
\begin{equation*}
\mathrm{clConv}_{\mathcal{F}}S=\{\omega\in\Omega\mid f(\omega)\leq \sup f(S)\text{ for all }f\in\mathcal{F}\}.
\end{equation*}
\end{definition}

\begin{lemma}\label{lem:convex}
Let $\mathcal{F}$ be the complete lattice cone generated by a linear space $\mathcal{A}$ of continuous functions on $\Omega$. Let $S\subset\Omega$. Then
\begin{equation*}
\mathrm{clConv}_{\mathcal{F}}S=\{\omega\in\Omega\mid a(\omega)\leq \sup a(S)\text{ for all }a\in\mathcal{A}\}.
\end{equation*}
\end{lemma}
\begin{proof}
Any element of $\mathcal{F}$ is a supremum of elements of $\mathcal{A}$, by Lemma \ref{lem:supgen}. The claim follows.
\end{proof}

\begin{proposition}\label{pro:convex}
Let $S\subset\Omega$. Then
\begin{equation*}
\mathrm{clConv}_{\mathcal{F}}S=\Phi^{-1}(\mathrm{clConv}_G\Phi(S)),
\end{equation*}
where $G$ is the complete lattice cone of convex, lower semi-continuous functions on $\mathcal{A}^*$ generated by evaluations on elements of $\mathcal{A}$.
\end{proposition}
\begin{proof}
By Remark \ref{rem:convex} and Lemma \ref{lem:convex} 
\begin{equation*}
\mathrm{clConv}_G\Phi(S)=\{a^*\in \mathcal{A}^*\mid a^*(a)\leq \sup \Phi(S)(a)\text{ for all }a\in \mathcal{A}\}.
\end{equation*}
Since for any $a\in \mathcal{A}$, $\Phi(S)(a)=a(S)$, we see that 
\begin{equation*}
\Phi^{-1}(\mathrm{clConv}_G\Phi(S))=\{\omega\in\Omega\mid a(\omega)\leq \sup a(S)\text{ for all }a\in \mathcal{A}\}.
\end{equation*}
An application of Lemma \ref{lem:convex} completes the proof.
\end{proof}

In view of Lemma \ref{lem:embed}, \ref{i:homeo}, and Proposition \ref{pro:convex} we define $\mathcal{F}$-convex hull of a set in the following  way.

\begin{definition}\label{def:conv}
Let $S\subset\Omega$. Let $\mathcal{A}$ be a linear space of continuous functions on $\Omega$ and let $\mathcal{F}$ be the complete lattice cone generated by $\mathcal{A}$. We define the $\mathcal{F}$-convex hull of $S$ by the  formula
\begin{equation*}
\mathrm{Conv}_{\mathcal{F}}S= \Phi^{-1}(\mathrm{Conv}_G\Phi(S)).
\end{equation*}
We shall say that $S\subset\Omega$ is $\mathcal{F}$-convex whenever $S=\mathrm{Conv}_{\mathcal{F}}S$.
Here $G$ is the complete lattice cone of convex, lower semi-continuous functions on $\mathcal{A}^*$ generated by evaluations on elements of $\mathcal{A}$.
\end{definition}

\begin{remark}
Note that above $\mathrm{Conv}_G\Phi(S)$ denotes the usual convex hull of a set $\Phi(S)\subset \mathcal{A}^*$. Moreover, if $S\subset\Omega$ is closed, then it is $\mathcal{F}$-convex if and only if $
S=\mathrm{clConv}_{\mathcal{F}}S$.
This is to say, Definition \ref{def:conv} and Definition \ref{def:convexhull} are consistent.
\end{remark}

\begin{lemma}\label{lem:sublevels}
Suppose that $f\in\mathcal{F}$. Then for any $t\in\mathbb{R}$ the sets $f^{-1}(-\infty,t))$ and $f^{-1}((-\infty,t])$ are $\mathcal{F}$-convex.
\end{lemma}
\begin{proof}
The claim follows by Defnition \ref{def:conv}, Lemma \ref{lem:embed}, \ref{i:fextend}  and an analogous claim for convex, lower semi-continuous functions on $\mathcal{A}^*$.
\end{proof}

\begin{proposition}\label{pro:convexequivalent}
Let $\mathcal{F}$ be a cone of continuous functions on $\Omega$ that contains constants. Let $K\subset\Omega$ be compact and let $C\geq 1$. Let $\omega\in \Omega$. The following conditions are equivalent:
\begin{enumerate}[label=(\roman*)]
\item\label{i:convoC} $f(\omega)\leq C \sup f(K)$ for all $f\in\mathcal{F}$,
\item\label{i:plus} $f(\omega)\leq C (\sup f(K))_+$ for all $f\in\mathcal{F}$,
\item\label{i:convo} $\omega\in\mathrm{clConv}_{\mathcal{F}}K$.
\end{enumerate} 
\end{proposition}
\begin{proof}
Clearly \ref{i:convoC} implies \ref{i:plus}. Suppose that \ref{i:plus} holds true.
Let $\mathcal{G}$ denote the complete lattice cone generated  by $\mathcal{F}$. Then also 
\begin{equation*}
g(\omega)\leq C( \sup g(K))_+\text{ for all }g\in\mathcal{G}.
\end{equation*}
Lemma \ref{lem:poussin} implies that if $\xi$ is a non-negative, convex and increasing function on $[0,\infty)$ then for any $f\in\mathcal{F}$, 
\begin{equation*}
\xi(f_+(\omega))\leq C\sup \xi(f_+(K)).
\end{equation*}
For $k\geq 2$, the function $t\mapsto t^k$ is convex and increasing on $[0,\infty)$. It follows that for $k\geq 2$
\begin{equation*}
f_+^k(\omega)\leq C \sup f_+^k(K).
\end{equation*}
Taking roots and letting $k$ tend to infinity yields
\begin{equation*}
f_+(\omega)\leq  \sup f_+(K).
\end{equation*}
Then applying the above to $f+\inf f(K\cup \{\omega\})\in\mathcal{F}$ shows that 
\begin{equation*}
f(\omega)\leq  \sup f(K).
\end{equation*}
That is, $\omega$ belongs to $\mathrm{clConv}_{\mathcal{F}}K$.
It is trivial that \ref{i:convo} implies \ref{i:convoC}.
\end{proof}

\begin{remark}\label{rem:nonequi}
Note that the above does not imply that for any  $K\subset\Omega$ and any $C\geq 1$ the set
$\mathrm{clConv}_{\mathcal{F}}K$ is equal to
\begin{equation*}
\{\omega\in\Omega\mid f(\omega)\leq C\sup\abs{f}(K)\text{  for all }f\in\mathcal{F}\}.
\end{equation*}
In particular, \ref{i:fconstantin} of Theorem \ref{thm:holoin} and \ref{i:foneRin} of Theorem \ref{thm:holoalgin} are not necessarily equivalent. 

Let $\mathcal{A}$ be the linear space of affine functions on an open, proper, convex set $\Omega$ of a finite-dimensional linear space. 
Let $K\subset\Omega$ be a compact, convex subset of $\Omega$, with non-empty interior. Then $K=\mathrm{clConv}_{\mathcal{A}}K$. However, the union  
\begin{equation*}
\bigcup\bigg\{\{\omega\in\Omega\mid a(\omega)\leq C\sup \abs{a}(K)\text{ for all }a\in\mathcal{A}\}\mid C\geq 1\bigg\}
\end{equation*}
consists of the elements of the affine hull of the compact set $K$. Since $K$ had non-empty interior, its affine hull contains $\Omega$. Thus, it is not equal to $K$.

In this example we see that \ref{i:foneRin} of Theorem \ref{thm:holoalgin} is satisfied, but \ref{i:fconstantin} of Theorem \ref{thm:holoin} is not.
\end{remark}

The proposition below shows that if $\mathcal{A}$ consists of real-parts of a complex algebra, then the situation described in Remark \ref{rem:nonequi} cannot take place.

\begin{proposition}\label{pro:fconvexity}
Let $\mathcal{A}$ consist of the real-parts of a complex algebra $\mathcal{B}$ of functions on $\Omega$, closed in $\tau_{\mathcal{K}}$ and containing constants. Then for any compact set $K\subset\Omega$ and any $C\geq 1$
\begin{equation*}
\mathrm{clConv}_{\mathcal{A}}K=\{\omega\in \Omega\mid a(\omega)\leq C\sup\abs{a}(K)\text{ for all }a\in\mathcal{A}\}.
\end{equation*}
Moreover, the above sets coincide with
\begin{equation*}
\{\omega\in \Omega\mid \abs{b(\omega)}\leq C\sup\abs{b}(K)\text{ for all }b\in\mathcal{B}\}.
\end{equation*}
\end{proposition}
\begin{proof}
Let 
\begin{equation*}
R=\{\omega\in \Omega\mid a(\omega)\leq C\sup\abs{a}(K)\text{ for all }a\in\mathcal{A}\} 
\end{equation*}
and 
\begin{equation*}
S= \{\omega\in \Omega\mid \abs{b(\omega)}\leq C\sup\abs{b}(K)\text{ for all }b\in\mathcal{B}\}.
\end{equation*}
Clearly, $\mathrm{clConv}_{\mathcal{A}}K\subset R$.

Let $\omega\in R$. Let $b\in\mathcal{B}$. Take $\theta\in\mathbb{C}$ of modulus one such that $\theta b(\omega)=\abs{b(\omega)}$. As $\omega\in R$ and $\theta b\in\mathcal{B}$ we see that 
\begin{equation*}
\abs{b(\omega)}\leq \sup C\abs{b}(K).
\end{equation*}
This shows that $R\subset S$.
Let now $\omega\in  S$. 
As $\mathcal{B}$ is an algebra, for any $b\in\mathcal{B}$ and $k\in\mathbb{N}$, $b^k$ belongs to $\mathcal{B}$. It follows that 
\begin{equation*}
\abs{b(\omega)}\leq C^{\frac1k} \sup \abs{b}(K).
\end{equation*}
Letting $k$ tend to infinity we get that 
\begin{equation*}
\abs{b(\omega)}\leq \sup \abs{b}(K).
\end{equation*}
Since $\mathcal{B}$ is closed in $\tau_{\mathcal{K}}$, for any $b\in\mathcal{B}$, $\exp(b)\in\mathcal{B}$. Applying the above to $\exp(b)$ and taking logarithms, we see that 
\begin{equation*}
\mathfrak{Re} b(\omega)\leq \sup \mathfrak{Re}b(K).
\end{equation*}
This shows that $\omega\in\mathrm{clConv}_{\mathcal{F}}K$ and concludes the proof.
\end{proof}

\begin{remark}
Let $\mathcal{A}$ be a linear space of continuous functions on $\Omega$  that contains constants. Then for any compact set $K\subset\Omega$, 
\begin{equation}\label{eqn:aconv}
\{\omega\in\Omega \mid \abs{a(\omega)}\leq \sup\abs{a}(K)\text{ for all }a\in\mathcal{A}\}=\mathrm{clConv}_{\mathcal{A}}K.
\end{equation}
Indeed, all functions in $\mathcal{A}$ are continuous, thus if $a\in\mathcal{A}$, then $a+\inf a(K)\in\mathcal{A}$. If $\omega$ belongs to the set on the left-hand side of (\ref{eqn:aconv}), then for all $a\in\mathcal{A}$, 
\begin{equation*}
a(\omega)+\inf  a(K)\leq \sup a(K)+\inf a(K),
\end{equation*}
so that $\omega$ belongs also to $\mathrm{clConv}_{\mathcal{A}}(K)$. The converse inclusion follows by Lemma \ref{lem:convex}.
\end{remark}

\section{Proofs}\label{s:proofs}

\subsection{Proof of Theorem \ref{thm:holoin}}

\begin{lemma}\label{lem:exhaustionC}
Let $\Omega$ be a topological space. Let $\mathcal{A}$ be a linear space of continuous functions on $\Omega$ that contains constants and let $\mathcal{F}$ be the complete lattice cone generated  by $\mathcal{A}$.
Suppose that $\Omega$ is exhaustible by compacta and that for any compact set $K\subset \Omega$ and $C\geq 1$, the set
\begin{equation*}
\{\omega\in\Omega\mid a(\omega)\leq C\sup\abs{a}(K)\text{ for all }a\in\mathcal{A}\}
\end{equation*}
is compact.
 Then there exists a continuous and non-negative $\mathcal{F}$-exhaustion function $p\colon\Omega\to\mathbb{R}$ such that
  \begin{equation*}
  p=\sup\{a_{\alpha}\mid \alpha\in A\}
 \end{equation*} 
 for some symmetric family $(a_{\alpha})_{\alpha\in A}$ of functions in $\mathcal{A}$. Moreover for any compact set $K$,  $p$ on $K$ is a maximum of a finite number of the functions in the family.
\end{lemma}
The lemma above is new, up to the author's knowledge, and the main new ingredient is the exploitation of the compactness of the sets  
\begin{equation*}
\{\omega\in\Omega\mid a(\omega)\leq C\sup\abs{a}(K)\text{ for all }a\in\mathcal{A}\}
\end{equation*}
for arbitrary $C\geq 1$, which in the setting of complex analysis is not widely studied.

\begin{proof}
We shall first choose a sequence of compacta $(K_i)_{i=1}^{\infty}$, such that
for each $i=1,2\dotsc,$
\begin{equation*}
\{\omega\in\Omega\mid \abs{a(\omega)}\leq 4^i\sup\abs{a}(K_i)\text{ for all }a\in\mathcal{A}\}\subset \mathrm{int}K_{i+1},
\end{equation*}
and such that
\begin{equation*}
\bigcup_{i=1}^{\infty}K_i=\Omega.
\end{equation*}
To this aim, $(L_i)_{i=1}^{\infty}$ be a sequence of compact subsets of $\Omega$ such that 
\begin{equation*}
\bigcup_{i=1}^{\infty}L_i=\Omega\text{ and }L_i\subset\mathrm{int}L_{i+1}\text{ for }i=1,2,\dotsc.
\end{equation*}
We shall construct the sets $(K_i)_{i=1}^{\infty}$ in such a way that for some strictly increasing function $j\colon\mathbb{N}\to\mathbb{N}$,  $L_{j(i-1)}\subset K_i\subset\mathrm{int}L_{j(i)}$ for $i=1,2,\dotsc$, where $j(0)=1$.
Suppose that $K_1,\dotsc,K_k$ are already constructed and $j(i)$ is defined for $i=0,\dotsc,k$. Set
\begin{equation*}
K_{k+1}=\{\omega\in\Omega\mid \abs{a(\omega)}\leq 4^k\sup\abs{a}( L_{j(k)})\text{ for all }a\in\mathcal{A}\}.
\end{equation*}
Then  $L_{j(k)}\subset K_{k+1}$ and $K_{k+1}$ is compact. Since it is covered by the interia of $(L_i)_{i=1}^{\infty}$, we may pick a finite subcover. Thus, it is contained in some set $L_{j(k+1)}$ of the family, and therefore 
\begin{equation*}
K_{k+1}\subset\mathrm{int}L_{j(k+1)}.
\end{equation*}
Now, since $j$ is strictly increasing, 
\begin{equation*}
\Omega=\bigcup_{k=1}^{\infty}L_{j(k)}=\bigcup_{k=1}^{\infty}K_k.
\end{equation*}
Moreover
\begin{equation*}
\{\omega\in\Omega\mid \abs{a(\omega)}\leq 4^i\sup\abs{a}(K_i)\text{ for all }a\in\mathcal{A}\}\subset \mathrm{int}L_{j(k+1)}\subset \mathrm{int}K_{i+1}.
\end{equation*}
Thus, the required sequence of the sets is constructed.

For each $i=1,2,\dotsc$ we shall now construct non-negative, continuous functions $p_i\in\mathcal{F}$ 
 such that 
\begin{equation*}
K_i\subset p_i^{-1}\big((-\infty, 1/2^i]\big)\text{ and }K_{i+2}\setminus K_{i+1}\subset p_i^{-1} \big([ 2^i,\infty)\big)\text{ for }i=1,2,\dotsc.
\end{equation*}
To this aim, observe that the compact set $K_{i+2}\setminus \mathrm{int}K_{i+1}$ is covered  by a family of open sets
\begin{equation*}
\{\omega\in \Omega\mid \abs{a(\omega)}> 4^i\sup\abs{a}(K_i)\}\text{ with }a\in\mathcal{A}.
\end{equation*}
By compactness we may pick a finite subcover and the  corresponding functions $a_1,\dotsc,a_{k_i}\in\mathcal{A}$. Let
\begin{equation*}
p_i=\max\Bigg\{\frac{\abs{a_l}}{2^i\norm{a_l}_{\mathcal{C}(K_i)}}\mid l=1,\dotsc,k_i\Bigg\}\vee 0.
\end{equation*}
Then  $p_i\in\mathcal{F}$ fulfils our requirements.
Now, let 
\begin{equation*}
p=\max\{p_i\mid i=1,2\dotsc\}.
\end{equation*}
Then $p$ belongs to $\mathcal{F}$. For $i=1,2,\dotsc$
\begin{equation*}
K_i\subset p_j^{-1}\big(\big(-\infty, 1/2^i\big]\big)\text{  for  all }j\geq i
\end{equation*}
and
\begin{equation*}
K_{i+2}\setminus \mathrm{int}K_{i+1}\subset p^{-1}\big([2^i,\infty)\big).
\end{equation*}
This shows that $p$ is proper. Moreover, we see that
\begin{equation*}
p=\max\{p_1,\dotsc,p_i\}\text{ on }K_i,
\end{equation*}
which immediately implies that $p$ is continuous.
Since each $p_i$ is a supremum of a symmetric family of functions in $\mathcal{A}$, thus so is $p$.  To  see that $p$ on any  compact is  a maximum of a finite family of functions in $\mathcal{A}$, it suffices to notice that each $p_i$ is such  a maximum, for $i=1,2,\dotsc$.

This completes the proof.
\end{proof}

\begin{proof}[Proof of Theorem \ref{thm:holoin}]
By Lemma \ref{lem:exhaustionC} and Remark \ref{rem:sigma} it  follows that \ref{i:fconstantin}  implies \ref{i:exhaustionfin}. Clearly, \ref{i:exhaustionfin} implies \ref{i:exhaustion} which in turn implies \ref{i:exhaustionbound}.

Let $\mathcal{F}$ denote the complete lattice cone generated by  $\mathcal{A}$.

Suppose that \ref{i:exhaustionbound} holds true. 
Let $p\colon\Omega\to [0,\infty)$ belong to $\mathcal{F}$. Then by Lemma \ref{lem:sublevels} for any $t\geq 0$ the set $p^{-1}([0,t])$ is $\mathcal{F}$-convex and compact. Suppose that $K\subset\Omega$ is compact. By the assumption, $p$ is bounded on compact sets, so there is some $t\geq 0$ such that $K\subset p^{-1}([0,t])$. Let $C\geq 1$. Consider the set 
\begin{equation*}
S=\{\omega\in\Omega\mid a(\omega)\leq C\sup\abs{a}(K)\text{ for all }a\in\mathcal{A}\}.
\end{equation*}
We claim that $S\subset p^{-1}([0,Ct])$. Indeed, by the assumption on $p\in\mathcal{F}$, there is a family $(a_{\alpha})_{\alpha\in A}$ of elements of $\mathcal{A}$ such that
\begin{equation*}
p=\sup\{\abs{ a_{\alpha}}\mid \alpha\in A\}.
\end{equation*}
Then $K\subset p^{-1}([0,t])$  implies that for $\alpha\in A$, $\sup \abs{a_{\alpha}}(K)\leq t$.
Let $\omega\in S$. Then
\begin{equation*}
p(\omega)=\sup \{\abs{a_{\alpha}}(\omega)\mid \alpha\in A\}\leq C\sup \{\sup\abs{a_{\alpha}}(K)\mid \alpha\in A\}\leq Ct.
\end{equation*}
This proves the claim and shows, due to properness of $p$, that \ref{i:fconstantin} is satisfied. 

Let us now show that \ref{i:exhaustionfin} implies existence of a family $(P_i)_{i=1}^{\infty}$ of compact, symmetric $\mathcal{A}$-polytopes that exhaust $\Omega$ and such that $P_i\subset \mathrm{int}P_{i+1}$ for $i=1,2,\dotsc$. 

By the assumption, $p$ on any compact set is a maximum of a finite family of functions in $\mathcal{A}$. Let 
\begin{equation*}
P_i=p^{-1}([0,i]).
\end{equation*}
These sets are compact, symmetric, $\mathcal{A}$-polytopes which exhaust $\Omega$. Moreover, as $p$ is non-negative and continuous
\begin{equation*}
 P_i=p^{-1}([0,i])\subset p^{-1}((-\infty, i+1))\subset \mathrm{int}P_{i+1}\text{ for }i=1,2,\dotsc.
 \end{equation*}
 This shows the existence of the requested $\mathcal{A}$-polytopes.
 
 If $\mathcal{A}$ consists of real-parts of a complex algebra, Proposition \ref{pro:fconvexity} shows that \ref{i:fconstantin} is equivalent to compactness of all sets $\mathrm{clConv}_{\mathcal{A}}K$, for compact $K\subset\Omega$. Since for any such $K$ there is some $i=1,2,\dotsc$ such that $K\subset P_i$, and these $P_i$ are $\mathcal{A}$-convex, we see that \ref{i:fconstantin} follows.
 
Under the assumption that $\Omega$ is separable, we shall show that \ref{i:fconstantin} is equivalent to $\Omega$ being complete with respect to $\mathrm{cl}\mathcal{A}$, where the closure is taken with respect to compact-open topology $\tau_{\mathcal{K}}$.

By Remark \ref{rem:metrisable} $\Omega$ is metrisable. Therefore, to prove that $\Omega$ is complete with respect to $\mathrm{cl}\mathcal{A}$, it suffices to consider Cauchy  sequences instead of general Cauchy nets.
Let us pick a Cauchy sequence $(\omega_n)_{n=1}^{\infty}$ of elements of $\Omega$.
The sequence induces a sequence $(\Phi(\omega_n))_{n=1}^{\infty}$ of linear, continuous functionals on $\mathrm{cl}\mathcal{A}$, which is pointwise bounded.
Therefore, by the Banach--Steinhaus uniform boundedness principle we can find a compact set $K\subset\Omega$ and a constant $C>0$ such that for all $a\in \mathcal{A}$  and all $n=1,2,\dotsc$
\begin{equation*}
a(\omega_{n})\leq C\norm{a}_{\mathcal{C}(K)}.
\end{equation*}
That is, the elements of  $(\omega_{n})_{n=1}^{\infty}$ belong to the set 
\begin{equation*}
\{\omega\in\Omega\mid a(\omega)\leq C\sup \abs{a}(K)\text{ for all }a\in\mathcal{A}\},
\end{equation*}
which is compact by \ref{i:fconstantin}. It follows that $(\omega_{n})_{n=1}^{\infty}$ converges. We have shown that $\Omega$ is complete with respect to $\mathrm{cl}\mathcal{A}$.

Suppose now conversely that $\Omega$ is complete with respect to $\mathrm{cl}\mathcal{A}$.
Recall that by Lemma \ref{lem:embed}, \ref{i:homeo}, the map $\Phi\colon\Omega\to\mathcal{A}^*$ is a homeomorphism onto its image, when $\mathcal{A}^*$ is equipped with the weak* topology. 
Let $K\subset\Omega$ and let $C\in\mathbb{R}$. Then the set 
\begin{equation*}
R=\{a^*\in\mathcal{A}^*\mid a^*(a)\leq C\sup \abs{a}(K)\text{ for all }a\in\mathcal{A}\}
\end{equation*}
is a compact subset of $\mathcal{A}^*$, by the Banach--Alaoglu theorem. 
Now 
\begin{equation*}
\{\omega\in\Omega\mid a(\omega)\leq C\sup \abs{a}(K)\text{ for all }a\in\mathcal{A}\}= \Phi^{-1}(R\cap \Phi(\Omega))
\end{equation*}
is therefore compact, as $\Phi(\Omega)$ is closed and $\Phi$ is a homeomorphism onto its image. This shows  that \ref{i:fconstantin} holds true.

Let us assume that \ref{i:fconstantin} is satisfied.  Suppose   that $\Omega$ is separable and that $\mathcal{A}$ consists of real-parts of a complex algebra. Remark \ref{rem:metrisable} shows that $\Omega$ is metrisable. Let $d$ be a metric on $\Omega$ that yields topology of  $\Omega$.
We shall show that $\Omega$ is an $(\mathrm{cl}\mathcal{A},d)$-space.

 Let $K_1,K_2,\dotsc$ be compact subsets of $\Omega$ such that
\begin{equation*}
\bigcup_{i=1}^{\infty}K_i=\Omega\text{ and }K_i\subset\mathrm{int}K_{i+1}\text{ for }i=1,2,\dotsc.
\end{equation*}
Taking their $\mathcal{F}$-convex hulls, and relabelling if necessary, we may moreover assume that these sets are $\mathcal{F}$-convex. We may moreover assume that none of these compact sets is equal to $\Omega$. Otherwise, the claim is trivial.

Let $\bar{\Omega}$ denote the completion of $\Omega$ with respect to $d$. Then $\bar{\Omega}$ is separable as a metrisable space with a separable dense subset $\Omega$. It admits a countable  basis $(U_i)_{i=1}^{\infty}$ of neighbourhoods of points in $\bar{\Omega}\setminus \Omega$. 
Let us take a function $k\colon\mathbb{N}\to\mathbb{N}$, with the property that each element of $\mathbb{N}$ appears infinitely often in its image.
As the sets $(K_j)_{j=1}^{\infty}$ are compact, for any $j=1,2,\dotsc$ there exists an element $\omega_{j}\in U_{k(j)}\cap \Omega\setminus K_j$. Relabelling the sets $(K_j)_{j=1}^{\infty}$ if necessary, we may assume that $\omega_j\in K_{j+1}$. 
By  Proposition \ref{pro:fconvexity} there exists $a_j\in\mathcal{A}$ such that
\begin{equation*}
\sup\abs{a_j}(K_j)\leq 2^{-j}\text{ and }a_j(\omega_j)>\sum_{i=1}^{j-1} \abs{a_i}(\omega_i)+j.
\end{equation*}
 Let us set
\begin{equation*}
a=\sum_{j=1}^{\infty}a_j.
\end{equation*}
Then the series converges in $\tau_{\mathcal{K}}$, so that $a\in\mathrm{cl}\mathcal{A}$.  Moreover, for each $j=1,2,\dotsc$
\begin{equation}\label{eqn:growth}
a(\omega_j)>a_j(\omega_j)-\sum_{i=1}^{j-1}\abs{a_i}(\omega_j)-\sum_{i=j+1}^{\infty}\abs{a_i}(\omega_j)\geq j-2^{-j}.
\end{equation}
We shall now show that $a$ verifies the fact that $\Omega$ is an $(\mathrm{cl}\mathcal{A},d)$-space. Let $(\omega_{\alpha})_{\alpha\in A}$ be a Cauchy net with respect to $d$ that does not converge to an element in $\Omega$, so that its limit $\omega\in\bar{\Omega}\setminus \Omega$.
By (\ref{eqn:growth}) we see that $\lim_{\alpha\in A}a(\omega_{\alpha})=\infty$, as $(a(\omega_{\alpha}))_{\alpha\in A}$ is unbounded on any neighbourhood of $\omega$. 

Suppose now that $\Omega$ is an $(\mathrm{cl}\mathcal{A},d)$-space for any metric $d$ that yields the  topology of $\Omega$. We shall show that $\Omega$ is complete with respect to $\mathrm{cl}\mathcal{A}$.  Remark \ref{rem:separable} shows that $\mathrm{cl}\mathcal{A}$ is separable in $\tau_\mathcal{K}$. Let $(a_i)_{i=1}^{\infty}$ be a dense and countable subset of $\mathrm{cl}\mathcal{A}$. Define for $\omega_1,\omega_2\in\Omega$
\begin{equation}\label{eqn:ddef}
    d(\omega_1,\omega_2)=\sum_{i=1}^{\infty}\frac1{2^i}\frac{\abs{a_i(\omega_1)-a_i(\omega_2)}}{\abs{a_i(\omega_1)-a_i(\omega_2)}+1}.
\end{equation}
Then $d\colon \Omega\times\Omega\to\mathbb{R}$ is a metric on $\Omega$. Since $\Omega$ is exhaustible by compacta, we may use an argument similar to the one in Lemma \ref{lem:embed} to show that above defined $d$ metrises the topology of $\Omega$. 

Pick now any Cauchy net $(\omega_{\alpha})_{\alpha\in A}$ with respect to $\mathrm{cl}\mathcal{A}$. 
 It is also a Cauchy net  with respect to $d$, so that we may take its equivalence class  $\omega\in\bar{\Omega}$. Suppose that $\omega$ is not equivalent to a point in $\Omega$.
  By the assumption, there exists $a\in\mathrm{cl}\mathcal{A}$ for which $\lim_{\alpha\in A}a(\omega_{\alpha})=\infty$.
 This is impossible, as the sequence is a Cauchy sequence with respect to $\mathrm{cl}\mathcal{A}$, so that $(a(\omega_{\alpha}))_{\alpha\in A}$ converges. Thus, $\Omega$ is complete with respect to $\mathrm{cl}\mathcal{A}$.
\end{proof}

\begin{remark}
Completeness of $\Omega$ implies that for any compact set $K\subset \Omega$ the set $\mathrm{clConv}_{\mathcal{F}}K$ is compact  in the following way.
Proposition \ref{pro:convex} tells us that for a compact set $K\subset\Omega$
\begin{equation*}
\mathrm{clConv}_{\mathcal{F}}K=\Phi^{-1}(\mathrm{clConv}_G\Phi(K)\cap\Phi(\Omega)).
\end{equation*}
The Mazur theorem and Lemma \ref{lem:embed} imply that $\mathrm{clConv}_G\Phi(K)$ is compact. 
The proof of Theorem \ref{thm:holoin} shows that if $\Omega$ is complete, then $\Phi(\Omega)$ is closed in the weak* topology on $\mathcal{A}^*$. Therefore
\begin{equation*}
\mathrm{clConv}_G\Phi(K)\cap\Phi(\Omega)
\end{equation*}
is a compact subset of $\Phi(\Omega)$.  As $\Phi$ is a homeomorphism onto its image, we see that $\mathrm{clConv}_{\mathcal{F}}K$ is compact subset of $\Omega$.
\end{remark}

\begin{remark}
Let us recall the following theorem of Krein and \v{S}mulian, \cite{Krein1940}.  
The original version was formulated for normed spaces. 
We shall be interested in a version for Fr\'echet spaces, see \cite[Theorem 3., p. 190, Definition 1., p. 184]{Wilansky1978}. We cite it below.
\end{remark}

\begin{theorem}\label{thm:krein-smulian}
Let $F$ be a Fr\'echet space. Then a convex set $Z\subset F^*$ is closed in the weak* topology if and only if its intersection with any polar set of an open neighbourhood in $F$ is weakly* compact. 
\end{theorem}

\begin{remark}
Let us note that equivalence of \ref{i:fconstantin} of Theorem \ref{thm:holoin} and of completeness of $\Omega$ with respect to $\mathcal{A}$ is analogous to the Krein--\v{S}mulian theorem.

If $\Phi(\Omega)\subset\mathcal{A}^*$ was convex, then the equivalence would follow from the Krein--\v{S}mulian theorem.
Indeed, in that case $\Phi(\Omega)$ is weakly* closed if and only if for any $C\geq 1$ and compact set $K\subset\Omega$ the set
\begin{equation*}
\{a^*\in\mathrm{cl}\mathcal{A}^*\mid a^*(a)\leq C\norm{a}_{\mathcal{C}(K)}\text{ for all }a\in\mathcal{A}\}\cap \Phi(\Omega)
\end{equation*}
is weakly* compact. This set is precisely the image, under $\Phi$, of the set 
\begin{equation*}
\{\omega\in\Omega\mid a(\omega)\leq C\sup\abs{a}(K)\text{ for all }a\in\mathcal{A}\},
\end{equation*}
which is compact by if \ref{i:fconstantin} of Theorem \ref{thm:holoin} holds true.
Completeness of $\Omega$ and weak* closedness of $\Phi(\Omega)$ are equivalent, by the definition of weak* topoloogy and by an application of the Banach--Steinhaus uniform boundedness principle.

The crucial observation in the proof of Theorem \ref{thm:holoin} is that $\Omega$ is metrisable, so that to prove its completeness it suffices to investigate behaviour of Cauchy sequences. Now, if $(\omega_n)_{n=1}^{\infty}$ is a Cauchy sequence, then $(\Phi(\omega_n))_{n=1}^{\infty}$ is a pointwise bounded sequence of linear functionals on $\mathcal{A}$. This would not be necessarily true, if we were to consider general nets of functionals on $\mathcal{A}$. 
Now, as $(\Phi(\omega_n))_{n=1}^{\infty}$ is pointwise bounded, we are free to employ the Banach--Steinhaus uniform boundedness principle and infer completeness. 
\end{remark}

If $\mathcal{A}$ is consists of real-parts of a closed complex  algebra, then there is an alternative proof of the fact that \ref{i:fconstantin} of Theorem \ref{thm:holoin} implies completeness of $\Omega$ with respect to $\mathrm{cl}\mathcal{A}$, that is analogous to the proof \cite[Theorem 7, p. 65]{Gunning1990}. It relies on the following lemma.

\begin{lemma}\label{lem:nonmetri}
Let $\Omega$ be a locally compact, $\sigma$-compact separable Hausdorff topological space. Suppose that $\mathcal{A}$ consists of real-parts of a complex algebra. Suppose that for any compact $K\subset \Omega$, $\mathrm{clConv}_{\mathcal{F}}K$ is compact. Then for any sequence $(\omega_n)_{n=1}^{\infty}$ with no accumulation point in $\Omega$ there exist a function $a\in\mathrm{cl}\mathcal{A}$ and a subsequence $(\omega_{n_j})_{j=1}^{\infty}$ of $(\omega_{n})_{n=1}^{\infty}$ such that  
\begin{equation*}
\lim_{j\to\infty}a(\omega_{n_j})=\infty.
\end{equation*}
\end{lemma}
\begin{proof}
Let $\mathcal{A}$ consist of real-parts of a complex algebra $\mathcal{B}$.

Let us pick a sequence  $(\omega_n)_{n=1}^{\infty}$ of elements of $\Omega$ with  no accumulation point in $\Omega$. Thus, for any compact $K\subset \Omega$, only finitely many elements of the sequence belong to $K$. 
We shall construct a sequence of compact sets $(K_i)_{i=1}^{\infty}$, for which 
\begin{equation*}
\bigcup_{i=1}^{\infty}K_i=\Omega\text{ and }K_i\subset \mathrm{int}K_{i+1}\text{ for }i=1,2,\dotsc.
\end{equation*}
Moreover, we shall find functions $(a_i)_{i=1}^{\infty}$ in $\mathcal{A}$ and a subsequence $(\omega_{n_j})_{j=1}^{\infty}$, such that  for  all $i=1,2,\dotsc$
\begin{equation*}
\omega_{n_i}\in  K_{i+1}\text{ and} \sup\abs{a_i}(K_i)\leq 2^{-i}
\end{equation*}
and
\begin{equation*}
a_i(\omega_{n_i})>i+\sum_{j=1}^{i-1}\abs{a_j}(\omega_{n_i}).
\end{equation*}
We shall proceed inductively.
Since $\Omega$ is exhaustible by compacta, there is a sequence of compact sets $(L_i)_{i=1}^{\infty}$ for which
\begin{equation}
\bigcup_{i=1}^{\infty}L_i=\Omega, L_i\subset \mathrm{int}L_{i+1}\text{ for }i=1,2,\dotsc.
\end{equation}
Suppose the functions $a_1,\dotsc,a_k\in\mathcal{A}$, elements $n_1,\dotsc,n_k$ and sets $K_1,\dotsc,K_k$ are chosen. There is an index $j_k\geq k$ such that $K_k\cup\{\omega_{n_1},\dotsc,\omega_{n_k}\}\subset\mathrm{int}L_{j_k}$.
By the assumption and Proposition \ref{pro:fconvexity}, the set
\begin{equation*}
K_{k+1}=\Big\{\omega\in\Omega\mid \abs{b}(\omega)\leq \sup \abs{b}(L_{j_k})\text{  for all }b\in\mathcal{B}\Big\}
\end{equation*}
is compact.
Since $(\omega_n)_{n=1}^{\infty}$ has no accumulation point, there is an element $\omega_{n_{k+1}}$ of the sequence that does not belong to $ K_{k+1}$.

Thus, there also exists $b_{k+1}\in\mathcal{A}$ for which $\abs{b_{k+1}}(\omega_{n_{k+1}})>\sup  \abs{b_{k+1}}(K_{k+1})$.
Normalising, taking power sufficiently high power and then real-part of a multiple of $b_{k+1}$ yields a function $a_{k+1}$ that we were after.  Since $j_k\geq k$, we see that 
\begin{equation*}
\bigcup_{i=1}^{\infty}K_i=\Omega.
\end{equation*}
We define now $a=\sum_{i=1}^{\infty}a_i$.
Since $\sup\abs{a_i}(K_j)\leq 2^{-i}$ for $i\geq j$, we see that the series converges in $\tau_{\mathcal{K}}$, and therefore $a\in\mathrm{cl}\mathcal{A}$.
However, for $j=1,2,\dotsc$, since $\omega_{n_j}\in K_i$ for $i\geq j+1$, so that $\abs{a_i}(\omega_{n_j})\leq 2^{-i}$, we have
\begin{equation*}
a(\omega_{n_j})\geq a_j(\omega_{n_j})-\sum_{i=1}^{j-1}\abs{a_i}(\omega_{n_j})-\sum_{i=j+1}^{\infty}\abs{a_i}(\omega_{n_j})\geq j-2^{-j}.
\end{equation*}
This shows  that $\lim_{j\to\infty} a(\omega_j)=\infty$.
\end{proof}

Remark \ref{rem:separable} and Remark \ref{rem:metrisable} show that $\Omega$ is metrisable. To prove that $\Omega$ is complete with respect to $\mathcal{A}$, it suffices to show that Cauchy sequences in $\Omega$ are convergent. There is no need to consider general Cauchy nets.  Lemma \ref{lem:nonmetri}  shows that any Cauchy sequence has an accumulation point, and therefore it is convergent in $\Omega$. Thus Lemma \ref{lem:nonmetri} proves that \ref{i:fconstantin} of  Theorem \ref{thm:holoin} implies completeness of $\Omega$. 

\subsection{Proof of Theorem \ref{thm:holoalgin}}

\begin{lemma}\label{lem:annulus}
Let $\Omega$ be a topological space. Let $\mathcal{F}$ be a cone of continuous functions on $\Omega$ that contains constants and satisfies the maximum principle. Let $K_1,K_2,K_3,K_4$ be compact and $\mathcal{F}$-convex subsets of $\Omega$ for which
\begin{equation*}
 K_1\subsetneq\mathrm{int}K_2, K_2\subset \mathrm{int}K_3\text{  and } K_3\subset K_4.
 \end{equation*}
Then there exist a non-negative, continuous function $p\colon\Omega\to\mathbb{R}$, that is a maximum of a finite family of functions in $\mathcal{F}$, such that
\begin{equation*}
K_1\subset  p^{-1}\big(\{0\}\big)
\end{equation*}
and 
\begin{equation*}
K_2\subset K_4\cap  p^{-1}\big((-\infty,1]\big)\subset\mathrm{int} K_3.
\end{equation*}
\end{lemma}
\begin{proof}
As  $K_1,K_2$ are compact, and $\mathcal{F}$-convex we see that for $i=1,2$
\begin{equation*}
K_i=\Big\{\omega\in\Omega\mid f(\omega)\leq\sup f(K_i)\text{ for all }f\in\mathcal{F}\Big\}.
\end{equation*}
Moreover, $K_4\setminus \mathrm{int}K_3$
is a compact set, covered by a collection of open sets
\begin{equation}
\Big\{\omega\in\Omega\mid f(\omega)> \sup f(K_2)\Big\},\text{ for all } f\in\mathcal{F},
\end{equation}
so it is likewise covered by a finite subcollection of such open sets. Let $f_1,\dotsc,f_k$ be the  corresponding functions in  $\mathcal{F}$.
Since  $\mathcal{F}$ satisfies the maximum principle and $K_1\subsetneq \mathrm{int}K_2$, we may assume that for $i=1,\dotsc,k$ we have
\begin{equation*}
\sup f_i(K_1)<\sup f_i(K_2).
\end{equation*}
 Then
\begin{equation}
\Big\{\omega\in K_4\mid \frac{f_i(\omega)-\sup f_i(K_1)}{\sup f_i(K_2)-\sup f_i(K_1)}\leq 1\text{ for }i=1,\dotsc,k\Big\}\subset\mathrm{int} K_3.
\end{equation}
It is readily verifiable that
\begin{equation}
p=\max\Big\{\frac{f_i-\sup f_i(K_1)}{\sup f_i(K_2)-\sup f_i(K_1)}\mid i=1,\dotsc,k\Big\}\vee 0
\end{equation}
satisfies our requirements.
\end{proof}

\begin{lemma}\label{lem:exhaustion}
Let $\Omega$ be a topological space. Let $\mathcal{F}$ be a cone of continuous functions on $\Omega$ that contains constants and satisfies the maximum principle. Let $\mathcal{G}$ be the  complete lattice cone generated by $\mathcal{F}$. Let $K_1,K_2,K_3,\dotsc$ be compact and $\mathcal{F}$-convex subsets of $\Omega$ for which
\begin{equation*}
 K_i\subsetneq\mathrm{int}K_{i+1}\text{ for }i=1,2,\dotsc\text{ and }\bigcup_{i=1}^{\infty}K_i=\Omega.
 \end{equation*}
 Then there exists a continuous, proper and non-negative function $p\in\mathcal{G}$ and  such that for any compact set $K$, $p$ on $K$ is a maximum of a finite number of functions in $\mathcal{F}$. 
\end{lemma}
\begin{proof}
For $i=1,2,\dotsc$ let $p_i\in\mathcal{G}$ be a non-negative, continuous function yielded by Lemma \ref{lem:annulus} and corresponding to the sets $K_i,K_{i+1},K_{i+2},K_{i+3}$.
Let 
\begin{equation}
p=\sup\{ip_i\mid i=1,2,\dotsc,\}.
\end{equation}
Observe that $p$ is non-negative and continuous. Indeed, if  $\omega\in K_j$ for $j=1,2,\dotsc$ then $p_i=0$ for $i\geq j$. Thus, $\Omega$ is covered by open sets, on each of which $p$ is a maximum of a finite number of non-negative and continuous functions.

We shall show that $p$ is proper. By the construction of Lemma \ref{lem:annulus}, $p_i$ is at least one on $K_{i+3}\setminus \mathrm{int}K_{i+2}$. Thus, for $i=1,2,\dotsc$
\begin{equation}
p\geq i\text{ on }\Omega\setminus\mathrm{int} K_{i+2}.
\end{equation}
This shows that for non-negative $t\in\mathbb{R}$ we have
\begin{equation*}
\{\omega\in\Omega\mid p(\omega)\leq t\}\subset \{\omega\in\Omega\mid p(\omega)\leq \lceil t\rceil\}\subset K_{\lceil t\rceil+3}.
\end{equation*}
By continuity of $p$ we see that preimages of compacta are compact, i.e., $p$ is proper.

This completes the proof.
\end{proof}

\begin{proof}[Proof of Theorem \ref{thm:holoalgin}]

Similarly to the argument in the proof of Theorem \ref{thm:holoin}, \ref{i:exhaustionboundR} implies \ref{i:foneRin}. 
Suppose now  that \ref{i:foneRin} is satisfied. If $\Omega$ is compact there is nothing to prove. Let us suppose that $\Omega$ is not compact.
Let $(K_i)_{i=1}^{\infty}$ be a sequence of compact, exhausting subsets of $\Omega$, cf. Remark \ref{rem:sigma}. By \ref{i:foneRin}, the sets $(\mathrm{clConv}_{\mathcal{F}}K_i)_{i=1}^{\infty}$ are compact and
\begin{equation*}
\bigcup_{i=1}^{\infty} \mathrm{clConv}_{\mathcal{F}}K_i=\Omega.
\end{equation*}
Since $\mathrm{int}K_i\subset\mathrm{intclConv}_{\mathcal{F}}K_i$  for $i=1,2,\dotsc$, we see that the sets $(\mathrm{clConv}_{\mathcal{F}}K_i)_{i=1}^{\infty}$ have non-empty interia that cover $\Omega$. Relabelling the sets if necessary, we may therefore assume that
\begin{equation*}
\mathrm{clConv}_{\mathcal{F}}K_i\subsetneq\mathrm{int}\mathrm{clConv}_{\mathcal{F}}K_{i+1}\text{ for }i=1,2\dotsc.
\end{equation*}
Indeed, if for some non-empty compact set $K$ we had $K=\mathrm{int}K$, then, due to assumption on connectedness, we would see that $K=\Omega$, contrary to the assumption that $\Omega$ is not compact.
Now, Lemma \ref{lem:exhaustion} shows that \ref{i:exhaustionfincon} holds true. Trivially,  this implies \ref{i:exhaustionconR}, which in turn implies \ref{i:exhaustionboundR}.

Let us now show that \ref{i:compactaf}  is equivalent to  \ref{i:foneRin}. Note that any $\mathcal{F}$-polytope is $\mathcal{F}$-convex by Lemma \ref{lem:sublevels}. If $K\subset\Omega$ is compact, then it is contained in some $\mathcal{F}$-polytope, and so does its $\mathcal{F}$-convex hull. Thus, \ref{i:foneRin} follows.  The converse implication follows thanks to an application of the equivalence of \ref{i:foneRin} and \ref{i:exhaustionfincon}. 
\end{proof}

\subsection{Proof of Theorem \ref{thm:levi}}
\begin{proof}[Proof of Theorem \ref{thm:levi}]
We shall employ Theorem \ref{thm:holoin} letting $\mathcal{A}$ denote the space of real-parts of functions in $\mathcal{B}$. Proposition \ref{pro:fconvexity} shows that holomorphic convexity of $\Omega$, i.e., compactness of $\mathrm{clConv}_{\mathcal{A}}K$ for any compact $K\subset\Omega$ is equivalent to \ref{i:fconstantin} of Theorem \ref{thm:holoin}.

Note that $\Omega$ is locally compact, $\sigma$-compact, separable and that $\mathcal{A}$ separates points of $\Omega$. 
Note also that $\mathcal{A}$ is closed in the compact-open topology, as $\mathcal{B}$ is. Thus all conditions of Theorem \ref{thm:holoin} are equivalent.

Note that $\Omega$ is complete with respect to $\mathcal{A}$ is and only if it is complete with respect to $\mathcal{B}$. Similarly, for any metric $d$, it is an $(\mathcal{A},d)$-space if and only if it is a $(\mathcal{B},d)$-space.

This shows that all conditions of Theorem \ref{thm:levi} are equivalent.
\end{proof}

\section{The lemma of Bremermann and Lelong implies the resolution of the Levi problem}\label{s:bremermann}

The lemma of Bremermann and of Lelong \cite[Theorem 2]{Bremermann1956} states that if $\Omega\subset\mathbb{C}^n$ is a domain of holomorphy, then any plurisubharmonic function on $\Omega$ is a Hartogs function. The set of Hartogs functions on $\Omega$, cf. \cite[Section 2]{Bremermann1956}, is the smallest convex cone of functions on $\Omega$ that contains the logarithms of moduli of holomorphic functions on $\Omega$ and is closed under taking:
\begin{enumerate}[label=(\roman*)]
    \item locally bounded suprema,
    \item monotonically decreasing limits.
\end{enumerate}
Moreover, the function belongs to the set of Hartogs functions on $\Omega$ if and only if its restrictions to any open, precompact subset $U$ of $\Omega$ is a Hartogs function on $U$, and any upper semi-continuous regularisation of a Hartogs function is again a Hartogs function, that is if $u\colon\Omega\to\mathbb{R}$ is a Hartogs function then also so is
\begin{equation*}
z\mapsto\limsup_{z'\to z}u(z').
\end{equation*}

\begin{remark}\label{rem:proof}
    We note that the proof of \cite[Theorem 2]{Bremermann1956} shows that any plurisubharmonic function $u\colon\Omega\to\mathbb{R}$ on a pseudoconvex domain $\Omega\subset\mathbb{C}^n$ can be written as
\begin{equation*}
    u(z)=\limsup_{z'\to z}\limsup_{m\to\infty} \frac{\log\abs{f_m(z')}}{m}
\end{equation*}
for some sequence $(f_m)_{n=1}^{\infty}$ of holomorphic functions on $\Omega$.
We refer also to \cite[Lemma 6.2]{Bremermann1959}.
\end{remark}

\begin{theorem}\label{thm:bremermann}
    The  Bremermann--Lelong lemma is equivalent to the positive resolution of the Levi problem.
\end{theorem}
\begin{proof}
    If the Levi problem has positive resolution, then \cite[Theorem 2]{Bremermann1956} shows that the lemma holds true. 
    
    Suppose that the converse is valid.   
 Let $\Omega\subset\mathbb{C}^n$ be pseudoconvex. We shall show that it is holomorphically convex. By Theorem \ref{thm:levi} we see that $\Omega$ admits a continuous plurisubharmonic exhaustion function $p\colon\Omega\to\mathbb{R}$, which is a Hartogs function thanks to our assumption. More specifically, see Remark \ref{rem:proof}, it can be written as 
 \begin{equation*}
    p(z)=\limsup_{z'\to z}\limsup_{m\to\infty} \frac{\log\abs{f_m(z')}}{m}
\end{equation*}
  for some sequence $(f_m)_{n=1}^{\infty}$ of holomorphic functions on $\Omega$.

Let $K\subset\Omega$ be compact.  Since $\Omega$ is open, there exists $\epsilon>0$ such that the set
\begin{equation*}
    K_{\epsilon}=\{\omega\in\Omega\mid \mathrm{dist}(\omega,K)\leq \epsilon\}
\end{equation*} 
is a compact subset of $\Omega$.
Since $p$ is continuous, it is bounded by $C$ on $K_{\epsilon}$. 
Therefore for any $z\in K$
\begin{equation*}
   \limsup_{m\to\infty} \frac{\log\abs{f_m(z)}}{m}\leq C.
\end{equation*}
Let $\mathcal{B}$ be the space of real-parts of holomorphic functions on $\Omega$ and let $\mathcal{S}$ be the convex cone generated by the logarithms of absolute values of holomorphic functions on $\Omega$.
The set
\begin{equation*}
   \mathrm{clConv}_{\mathcal{B}}K= \{\omega\in\Omega\mid a(\omega)\leq \sup a(K)\text{ for all }a\in\mathcal{B}\}
\end{equation*}
is contained in the set 
\begin{equation*}
    L=\{\omega\in\Omega\mid b(\omega)\leq \sup b(K)\text{ for all }b\in\mathcal{S}\}.
\end{equation*}
If $\omega\in L$, i.e., $b(\omega)\leq \sup b(K)$ for all $b\in\mathcal{S}$, then  by the Hahn--Banach theorem for convex cones, \cite[Theorem 2.1]{Roth2000}, and by the Riesz' representation theorem of non-negative functionals there exists a probability measure $\mu_{\omega}$ on $K$ such that for all $b\in\mathcal{S}$
\begin{equation*}
    b(\omega)\leq \int_Kb\, d\mu_{\omega}.
\end{equation*}
Therefore for all $m=1,2,\dotsc$ we have
\begin{equation*}
    \frac1m\log\abs{f_m}(\omega)\leq \int_K\frac1m\log\abs{f_m}\, d\mu_{\omega}.
\end{equation*}
Thus
\begin{equation*}
    \limsup_{m\to\infty} \frac1m\log\abs{f_m}(\omega)\leq \limsup_{m\to\infty}\Bigg(\int_K\ \frac1m\log\abs{f_m}(\omega) \, d\mu_{\omega}\Bigg).
\end{equation*}
Now,  thanks to the boundedness of $p$ and thus also of $\limsup_{n\to\infty}\frac1m\log\abs{f_m}$ on $K$ we can employ the Fatou lemma and see that 
\begin{equation*}
   p(\omega)= \limsup_{m\to\infty} \frac1m\log\abs{f_m}(\omega)\leq \Bigg(\int_K\ \limsup_{m\to\infty}\frac1m\log\abs{f_m}(\omega) \, d\mu_{\omega}\Bigg)\leq C.
\end{equation*}
This completes the proof, as $p$ has compact sublevel sets and therefore $\mathrm{clConv}_{\mathcal{B}}K$ is compact for any compact $K\subset\Omega$. This is to say, $\Omega$ is holomorphically convex.
\end{proof}


\end{document}